\documentclass{amsart}
\usepackage{xy}
\xyoption{all}
\usepackage{graphicx}
\usepackage{amsmath}
\usepackage{amsfonts}
\newtheorem{theorem}{Theorem}[section]

\newtheorem{corollary}[theorem]{Corollary}

\newtheorem{lemma}[theorem]{Lemma}
\newtheorem{proposition}[theorem]{Proposition}
\theoremstyle{remark}
\newtheorem{remark}[theorem]{Remark}
\numberwithin{equation}{section}

\begin{document}
\title[The Geometry of $3$-quasi-Sasakian manifolds]{The Geometry of $3$-quasi-Sasakian manifolds}

\author[B. Cappelletti Montano]{Beniamino Cappelletti Montano}
 \address{Dipartimento di Matematica,
Universit\`{a} degli Studi di Bari, Via E. Orabona 4, 70125 Bari, Italy}
 \email{cappelletti@dm.uniba.it}

\author[A. De Nicola]{Antonio De Nicola}
 \address{CMUC, Department of Mathematics, University of Coimbra, 3001-454 Coimbra, Portugal}
 \email{antondenicola@gmail.com}

\author[G. Dileo]{Giulia Dileo}
 \address{Dipartimento di Matematica,
Universit\`{a} degli Studi di Bari, Via E. Orabona 4, 70125 Bari, Italy}
 \email{dileo@dm.uniba.it}



\begin{abstract}
$3$-quasi-Sasakian manifolds were studied systematically by the
authors in a recent paper as a suitable setting unifying
$3$-Sasakian and $3$-cosymplectic geometries. This paper throws new
light on their geometric\allowbreak structure which appears to be generally
richer compared to the $3$-Sasakian subclass. In fact, it turns out
that they are multiply foliated by four distinct fundamental
foliations. The study of the transversal geometries with respect to
these foliations allows us to link the $3$-quasi-Sasakian manifolds
to the more famous hyper-K\"{a}hler and quaternionic-K\"{a}hler
geometries. Furthermore, we strongly improve the splitting results
previously obtained; we prove that any $3$-quasi-Sasakian manifold of
rank $4l+1$ is $3$-cosymplectic and any $3$-quasi-Sasakian manifold
of maximal rank is $3$-$\alpha$-Sasakian.
\end{abstract}

\maketitle

\section{Introduction}
The well-known classes of $3$-Sasakian and $3$-cosymplectic
manifolds belong to the wider family of almost $3$-contact metric
manifolds. Nevertheless, both classes sit also perfectly into the
narrower class of $3$-quasi-Sasakian manifolds which, as we will
see, is a very natural framework for a unified study of the
aforementioned geometries. A similar chain of inclusions takes place
in the case of a single almost contact metric structure, whereas the
class of quasi-Sasakian manifolds encloses both Sasakian and
cosymplectic manifolds, but in the setting of $3$-structures the
interrelations between the triples of tensors produce key additional
properties making the choice of the $3$-quasi-Sasakian framework
still more natural. $3$-quasi-Sasakian manifolds were introduced
long ago but their first systematic study was carried out by the
authors in \cite{mag}. There, it was proven that in any
$3$-quasi-Sasakian manifold
$(M,\phi_\alpha,\xi_\alpha,\eta_\alpha,g)$ of dimension $4n+3$ the
vertical distribution $\mathcal V$ generated by the three Reeb
vector fields is completely integrable determining a canonical
totally geodesic and Riemannian foliation. The characteristic vector
fields obey the commutation relations
$[\xi_\alpha,\xi_\beta]=c\xi_\gamma$ for any even permutation
$\left(\alpha,\beta,\gamma\right)$ of $\left\{1,2,3\right\}$ and
some $c \in \mathbb R$. Furthermore, it was shown that the ranks of
the $1$-forms $\eta_1,\eta_2,\eta_3$ coincide giving a single
well-defined rank which falls into one of two possible families:
$4l+3$ or $4l+1$ for some $0 \leq l \leq n$. A splitting theorem was
proven for the manifolds in the first of the two families just for
the case $c=2$, under some additional hypotheses, while a sufficient
condition for those of rank $4l+1$ to be $3$-cosymplectic was found.
In this paper, the whole geometric structure of $3$-quasi-Sasakian
manifolds is enlightened. They appear to have in general an even
richer structure compared to that of the more famous $3$-Sasakian
subclass. In fact, it turns out that they are multiply foliated by
four distinct fundamental foliations, three of which become trivial
in the case of $3$-Sasakian manifolds. The study of the transversal
geometries with respect to these foliations allows us to link
$3$-quasi-Sasakian manifolds to the more famous hyper-K\"{a}hler and
quaternionic-K\"{a}hler geometries. Furthermore, we strongly
improve the splitting results previously found, showing that any
$3$-quasi-Sasakian manifold of rank $4l+3$ is locally the Riemannian
product of a $3$-$\alpha$-Sasakian manifold and a hyper-K\"{a}hler
manifold without any additional hypothesis. We obtain many
additional properties of $3$-quasi-Sasakian manifolds,
characterizing those of minimal and maximal rank and we study
Riemannian and Ricci curvature,  determining exactly which
$3$-quasi-Sasakian manifolds are Einstein or $\eta$-Einstein. Some
topological obstructions to the existence of $3$-quasi-Sasakian
structures on a given manifold are also found.

The article is organized as follows. In $\S$\!\! $2$ we  briefly recall
the required preliminaries about almost contact metric geometry and
$3$-structures, the two pillars supporting
$3$-quasi-Sasakian geometry. The most relevant results already known
about $3$-quasi-Sasakian manifolds are also summarized. In the third
section we mainly prove that all $3$-quasi-Sasakian manifolds of
rank $4l+1$ are $3$-cosymplectic. It follows that a
$3$-quasi-Sasakian manifold is Ricci-flat if and only if it is
$3$-cosymplectic. Such a corollary may be thought as an
odd-dimensional analogue of the well-known fact that any
quaternionic-K\"{a}hler manifold is Ricci-flat if and only if it is
locally hyper-K\"{a}hler. $\S$\! $4$ is devoted to the study of the
complementary class: $3$-quasi-Sasakian manifolds of rank $4l+3$. We
show that any $3$-quasi-Sasakian manifold of maximal rank is
$3$-$\alpha$-Sasakian (cf. \cite{friedrich}) and it is $3$-Sasakian
if and only if the constant $c$ in the commutators
$[\xi_\alpha,\xi_\beta]=c\xi_\gamma$ is equal to $2$. Some results
concerning the Riemannian and sectional curvatures of
$3$-quasi-Sasakian manifolds of rank $4l+3$ are also obtained.
Finally, $\S$\! $5$ contains the main results concerning the geometry
of $3$-quasi-Sasakian manifolds, shaped by the four fundamental
foliations which we find out to be canonically associated to each of
them. We start by analyzing the vertical foliation $\mathcal V$. The
use of an adapted connection (cf. \cite{cappellettidenicola})
derived from the Bott connection allows us to show that the
canonical transversal structure with respect to $\mathcal V$ is
projectable (as  a whole) and it is almost quaternionic-Hermitian.
The projectability of each structure tensor $\phi_\alpha$ with
respect to $\mathcal V$ is then shown to be equivalent to the
integrability of the horizontal distribution. This happens if and
only if the $3$-quasi-Sasakian manifold $M$ has minimal rank $3$ and
in that case we prove that $M$ is locally the Riemannian product of
a hyper-K\"{a}hler manifold and a three dimensional sphere. Next, we
study a second fundamental integrable distribution, denoted by
${\mathcal E}^{4m}$ (where $m=n-l$), in any $3$-quasi-Sasakian
manifold of rank $4l+3$. The leaves turn out to be hyper-K\"{a}hler
while the leaf space is $3$-$\alpha$-Sasakian. Then, we prove that
the integrability of the distribution ${\mathcal E}^{4m+3}={\mathcal
E}^{4m}\oplus {\mathcal V}$ gives rise to a foliation of $M$ in
$3$-quasi-Sasakian leaves of minimal rank whose leaf space has
quaternionic-K\"{a}hler structure. In this way we generalize to the
class of $3$-quasi-Sasakian manifold a fundamental result proven by
Ishihara \cite{ishihara} for $3$-Sasakian manifolds with respect to
the foliation ${\mathcal V}$ which turned out to be fundamental for
the subsequent studies of C. P. Boyer, K. Galicki and many others
giving to that class of manifolds their current relevance. Finally,
in Theorem \ref{profinalissimo} we show that every
$3$-quasi-Sasakian manifold of rank $4l+3$ admits a canonical
transversal hyper-K\"{a}hler structure given by a foliation whose
leaves are $3$-$\alpha$-Sasakian. As a very important corollary we
are able to greatly improve the known splitting result obtained in
\cite{mag}. We prove that any $3$-quasi-Sasakian manifold of rank
$4l+3$ is locally the Riemannian product of a $3$-$\alpha$-Sasakian
manifold and a hyper-K\"{a}hler manifold. We do not make any
assumption neither on the metric nor on the constant $c$. A number
of notable consequences follow, among which the reduction of the
structure group to $\emph{Sp}(m)\times\emph{Sp}(l) \times I_3$ and
the existence of nine contact-symplectic pairs on any
$3$-quasi-Sasakian manifold of rank $4l+3$. It also allows us to
compute the full Ricci tensor, proving that any $3$-quasi-Sasakian
manifold of rank $4l+3$ has positive scalar curvature and permits to
determine which $3$-quasi-Sasakian manifolds are Einstein or
$\eta$-Einstein.

\section{Preliminaries}
An \emph{almost contact manifold} is an odd-dimensional manifold $M$
which carries a tensor field $\phi$ of type $(1,1)$, a vector field
$\xi$, called \emph{characteristic} or \emph{Reeb vector field}, and
a $1$-form $\eta$ satisfying $\phi^2=-I+\eta\otimes\xi$ and
$\eta\left(\xi\right)=1$, where $I \colon TM\rightarrow TM$ is the
identity mapping. From the definition it follows also that
$\phi\xi=0$, $\eta\circ\phi=0$ and that the $(1,1)$-tensor field
$\phi$ has constant rank $2n$ (cf. \cite{blair1}). An almost contact
manifold is said to be \emph{normal} if the tensor field
$N^{(1)}=[\phi,\phi]+2d\eta\otimes\xi$ vanishes identically. It is
known that any almost contact manifold
$\left(M,\phi,\xi,\eta\right)$ admits a Riemannian metric $g$ such
that
$g\left(\phi\,\cdot,\phi\,\cdot\right)=g\left(\cdot,\cdot\right)-\eta\otimes\eta$
holds. This metric, in general not unique, is called a
\emph{compatible metric} and the manifold $M$ together with the
structure $\left(\phi,\xi,\eta,g\right)$ is called an \emph{almost
contact metric manifold}. As an immediate consequence one has
$\eta=g\left(\xi,\cdot\right)$. The $2$-form $\Phi$ on $M$ defined
by $\Phi\left(X,Y\right)=g\left(X,\phi Y\right)$ is called the
\emph{fundamental $2$-form} of the almost contact metric manifold
$M$.
The following formula gives the expression of the covariant
derivative of $\phi$ in terms of the remaining structure tensors in
any almost contact metric manifold (\cite{blair1}) and it will be
useful in the sequel,
\begin{equation}\label{olszakformula}
\begin{split}
2g((\nabla_X\phi)Y,Z)&=3d\Phi(X,\phi Y,\phi
Z)-3d\Phi(X,Y,Z)+g(N^{(1)}(Y,Z),\phi X)\\
&\quad+N^{(2)}(Y,Z)\eta(X)+2d\eta(\phi Y,X)\eta(Z)-2d\eta(\phi
Z,X)\eta(Y),
\end{split}
\end{equation}
where $N^{(2)}$ is the tensor defined by $N^{(2)}(X,Y)=({\mathcal L}_{\phi X}\eta)(Y)-({\mathcal L}_{\phi Y}\eta)(X)$ (cf. \cite{blair1}).

Almost contact metric manifolds such that both $\eta$ and $\Phi$ are
closed are called \emph{almost cosymplectic manifolds} and almost
contact metric manifolds such that $d\eta=\Phi$ are called
\emph{contact metric manifolds}. Finally, a normal almost
cosymplectic manifold is called a \emph{cosymplectic manifold} and a
normal contact metric manifold is said to be a \emph{Sasakian
manifold}.

The notion of  quasi-Sasakian structure, introduced by D. E. Blair
in \cite{blair0}, unifies those of Sasakian and cosymplectic
structures. A \emph{quasi-Sasakian manifold} is a normal almost
contact metric manifold such that $d\Phi=0$. A quasi-Sasakian
manifold $M$ (or more generally an almost contact manifold) of
dimension $2n+1$ is said to be of rank $2p$ (for some $p\leq n$)  if
$\left(d\eta\right)^p\neq 0$ and $\eta\wedge\left(d\eta\right)^p=0$
on $M$, and to be of rank $2p+1$ if
$\eta\wedge\left(d\eta\right)^p\neq 0$ and
$\left(d\eta\right)^{p+1}=0$  on $M$ (cf. \cite{blair0,tanno}). It
was proven in \cite{blair0} that there are no quasi-Sasakian
manifolds of even rank. Let the rank of $M$ be $2p+1$. Then, the
tangent bundle of $M$ splits into two subbundles as follows:
$TM={\mathcal E}^{2p+1}\oplus{\mathcal E}^{2q}$, $p+q=n$, where
\begin{equation*}
{\mathcal E}^{2q}=\{X\in TM\; | \;  i_X \eta=0  \mbox{ and } i_X
d\eta=0\}
\end{equation*}
and  ${\mathcal E}^{2p+1}={\mathcal
E}^{2p}\oplus\left\langle\xi\right\rangle$, ${\mathcal E}^{2p}$
being the orthogonal complement of ${\mathcal
E}^{2q}\oplus\left\langle\xi\right\rangle$ in $TM$. These
distributions satisfy $\phi {\mathcal E}^{2p}={\mathcal E}^{2p}$ and
$\phi {\mathcal E}^{2q}={\mathcal E}^{2q}$ (cf. \cite{tanno}).
Notice that for all $x\in M$ the subspace ${\mathcal E_x}^{2q}$
 coincides with the characteristic system defined by E. Cartan in \cite{cartan}
for an arbitrary differential form. The class of a differential form
is one of the integral invariants defined by Cartan. The codimension
$2p+1$ of ${\mathcal E_x}^{2q}$ is called by Cartan the
\textit{class} of $\eta$ in $x$. It is easy to verify that when the
class of $\eta$ is constant the characteristic system has constant
rank in any point and the determined distribution is integrable.
This is the case in all important examples of quasi-Sasakian
manifolds, such as Sasakian and cosymplectic manifolds. Thus, we
will only consider, as Blair and Tanno  implicitly did,
quasi-Sasakian manifolds of constant class, i.e. of fixed (odd)
rank. So, the rank of Blair and Tanno coincides with the class of
Cartan.

Some useful properties of quasi-Sasakian manifolds will be now mentioned.
For a quasi-Sasakian manifold we have the relation (cf. \cite{olszak2})
\begin{equation}\label{formulaquasisasaki}
\left(\nabla_{X}\phi\right)Y=-g\left(\nabla_{X}\xi,\phi
Y\right)\xi-\eta\left(Y\right)\phi\nabla_{X}\xi,
\end{equation}
which generalizes the well-known conditions $\nabla\phi=0$ and
$\left(\nabla_X\phi\right)Y=g\left(X,Y\right)\xi-\eta\left(Y\right)X$
characterizing respectively cosymplectic and Sasakian manifolds.
The quasi-Sasakian condition reflects also in some properties of
curvature and of the Reeb vector field. In fact we have the
following results.

\begin{lemma}[\cite{blair0},\cite{olszak2}]
Let $\left(M,\phi,\xi,\eta,g\right)$ be a quasi-Sasakian manifold.
Then
\begin{enumerate}
  \item[(i)] the Reeb vector field $\xi$ is Killing and its integral curves are geodesics;
  \item[(ii)] the Ricci curvature in the direction of $\xi$ is given by
\begin{equation}\label{olszakricci}
\emph{Ric}\left(\xi\right)=\|\nabla\xi\|^2.
\end{equation}
\end{enumerate}
\end{lemma}

We now come to the main topic of our paper, i.e. $3$-quasi-Sasakian
geometry, which is framed into the more general setting of almost
$3$-contact geometry. An  \emph{almost $3$-contact manifold}  is a
$\left(4n+3\right)$-dimensional smooth  manifold $M$ endowed with
three almost contact structures $\left(\phi_1,\xi_1,\eta_1\right)$,
$\left(\phi_2,\xi_2,\eta_2\right)$,
$\left(\phi_3,\xi_3,\eta_3\right)$ satisfying the following
relations, for any even permutation
$\left(\alpha,\beta,\gamma\right)$ of $\left\{1,2,3\right\}$,
\begin{equation}
\begin{split}\label{3-sasaki}
\phi_\gamma=\phi_{\alpha}\phi_{\beta}-\eta_{\beta}\otimes\xi_{\alpha}=-\phi_{\beta}\phi_{\alpha}+\eta_{\alpha}\otimes\xi_{\beta},\quad\\
\xi_{\gamma}=\phi_{\alpha}\xi_{\beta}=-\phi_{\beta}\xi_{\alpha}, \quad
\eta_{\gamma}=\eta_{\alpha}\circ\phi_{\beta}=-\eta_{\beta}\circ\phi_{\alpha}.
\end{split}
\end{equation}
This notion was introduced by Y. Y. Kuo (\cite{kuo}) and,
independently, by C. Udriste (\cite{udriste}). In \cite{kuo} Kuo
proved that given an almost contact $3$-structure
$\left(\phi_\alpha,\xi_\alpha,\eta_\alpha\right)$, there exists a
Riemannian metric $g$ compatible with each of them and hence we can
speak of \emph{almost contact metric $3$-structures}. It is well
known that in any almost $3$-contact metric manifold the Reeb vector
fields $\xi_1,\xi_2,\xi_3$ are orthonormal with respect to the
compatible metric $g$ and that the structural group of the tangent
bundle is reducible to $\emph{Sp}\left(n\right)\times I_3$.
Moreover, by putting
${\mathcal{H}}=\bigcap_{\alpha=1}^{3}\ker\left(\eta_\alpha\right)$
one obtains a $4n$-dimensional distribution on $M$ and the tangent
bundle splits as the orthogonal sum
$TM={\mathcal{H}}\oplus{\mathcal{V}}$, where ${\mathcal
V}=\left\langle\xi_1,\xi_2,\xi_3\right\rangle$. We will call any
vector belonging to the distribution $\mathcal H$ \emph{horizontal}
and any vector belonging to the distribution $\mathcal V$
\emph{vertical}. An almost $3$-contact manifold $M$ is said to  be
\emph{hyper-normal} if each  almost contact structure
$\left(\phi_\alpha,\xi_\alpha,\eta_\alpha\right)$ is
normal. 

A \emph{$3$-quasi-Sasakian manifold} is, by definition, an almost
$3$-contact metric manifold such that each structure
$(\phi_\alpha,\xi_\alpha,\eta_\alpha,g)$ is quasi-Sasakian.
Important subclasses of the above defined class are the well-known
$3$-Sasakian and $3$-cosymplectic manifolds. Many results about
$3$-quasi-Sasakian manifolds have been found in \cite{mag}.

\begin{theorem}[\cite{mag}]\label{principale}
Let $(M,\phi_\alpha,\xi_\alpha,\eta_\alpha,g)$ be a
$3$-quasi-Sasakian manifold. Then the distribution spanned by the
Reeb vector fields $\xi_1$, $\xi_2$, $\xi_3$ is integrable and
defines a totally geodesic and Riemannian foliation $\mathcal V$ of
$M$. In particular, we have, for an even permutation
$(\alpha,\beta,\gamma)$ of $\left\{1,2,3\right\}$, that
$[\xi_\alpha,\xi_\beta]=c\xi_\gamma$ for some $c\in\mathbb{R}$.
\end{theorem}

According to Theorem \ref{principale}, the geometry of
$3$-quasi-Sasakian manifolds with $c=0$ and
those with $c\neq 0$, is very different.
This can be seen, for instance, in the notion of the ``rank'' of a
$3$-quasi-Sasakian manifold, which is well defined due to the following
theorem.

\begin{theorem}[\cite{mag}]\label{principale1}
Let $(M,\phi_\alpha,\xi_\alpha,\eta_\alpha,g)$ be a
$3$-quasi-Sasakian manifold of dimension $4n+3$. Then the $1$-forms
$\eta_1$, $\eta_2$, $\eta_3$ have the same rank, which is called the
\emph{rank} of the $3$-quasi-Sasakian manifold $M$. Furthermore,
this rank is equal to $4l+1$ or $4l+3$, for some $l\leq n$,
according to $c=0$ or $c\ne0$ respectively.
\end{theorem}

Now we collect some results on $3$-quasi-Sasakian manifolds, which
we will use in the sequel. As before, we refer the reader to \cite{mag}
for the details.

\begin{proposition}\label{lemmi}
In any $3$-quasi-Sasakian manifold
$(M,\phi_\alpha,\xi_\alpha,\eta_\alpha,g)$ we have:
\begin{enumerate}
  \item[(i)] $d\eta_\alpha(X,\xi_\beta)=0$ for all $X\in\Gamma\left(\mathcal
  H\right)$ and $\alpha,\beta\in\left\{1,2,3\right\}$;
  \item[(ii)] every Reeb vector field $\xi_\alpha$ is an
  infinitesimal automorphism with respect to the distribution $\mathcal
  H$;
  \item[(iii)] $d\eta_\alpha=\frac{1}{2}{\mathcal
  L}_{\xi_\beta}\Phi_\gamma$, for any even permutation
  $(\alpha,\beta,\gamma)$ of $\left\{1,2,3\right\}$;
  \item[(iv)] $d\eta_\alpha\left(X,\phi_\alpha
Y\right)=d\eta_\beta\left(X,\phi_\beta Y\right)$ for all
$X,Y\in\Gamma\left(\mathcal H\right)$ and
$\alpha,\beta\in\left\{1,2,3\right\}$;
  \item[(v)] $d\eta_\alpha\left(\phi_\beta X,\phi_\beta
Y\right)=-d\eta_\alpha\left(X,Y\right)$  for all
$X,Y\in\Gamma\left(\mathcal H\right)$ and $\alpha\neq\beta$;
\item[(vi)] $d\eta_\alpha(\phi_\beta X,Y)=d\eta_\gamma(X,Y)$ for all
$X,Y\in\Gamma\left(\mathcal H\right)$ and for any even permutation
  $(\alpha,\beta,\gamma)$ of $\left\{1,2,3\right\}$.
\end{enumerate}
\end{proposition}

Finally, for the Levi Civita connection of a $3$-quasi-Sasakian
manifold, we have $\nabla_{\xi\alpha}\xi_\alpha=0$ and
$\nabla_{\xi\alpha}\xi_\beta=\frac{c}{2}\xi_\gamma$ for any even
permutation $(\alpha,\beta,\gamma)$ of $\left\{1,2,3\right\}$.
Hence, the Riemannian curvature satisfies
\begin{equation}\label{Rxi}
R(\xi_\alpha,\xi_\beta)\xi_\beta=\frac{c^2}{4}\xi_\alpha\quad\mbox{and}\quad R(\xi_\alpha,\xi_\beta)\xi_\gamma=0
\end{equation}
for an even permutation
$(\alpha,\beta,\gamma)$ of $\left\{1,2,3\right\}$.

\section{Further results on $3$-quasi-Sasakian manifolds}
In  the  following \ we  will  use \  the  notation \ ${\mathcal
E}^{4m}:=\{X\in{\mathcal H}\; | \;  i_X d\eta_\alpha=0 \   \textrm{
for some } \alpha\in\{1,2,3\}  \}$, while ${\mathcal E}^{4l}$ will
be the orthogonal complement of ${\mathcal E}^{4m}$ in ${\mathcal
H}$, ${\mathcal E}^{4l+3}:={\mathcal E}^{4l} \oplus {\mathcal V}$,
and ${\mathcal E}^{4m+3}:={\mathcal E}^{4m} \oplus{\mathcal V}$. It
is easy to see that $\phi_\alpha({\mathcal E}^{4m})={\mathcal
E}^{4m}$ and $\phi_\alpha({\mathcal E}^{4l})={\mathcal E}^{4l}$ for
each $\alpha\in\left\{1,2,3\right\}$. Note also that, regarding to
the definition of ${\mathcal E}^{4m}$, if for a horizontal vector
$X$ $i_{X}d\eta_\alpha=0$ for \emph{some}
$\alpha\in\left\{1,2,3\right\}$, then by Lemma 5.4 in \cite{mag}
$i_{X}d\eta_\delta=0$ for \emph{any}
$\delta\in\left\{1,2,3\right\}$. We remark that in the case of
$3$-quasi-Sasakian manifolds of rank $4l+3$ the distribution
${\mathcal E}^{4m}$ is integrable since it coincides with the
distribution defined by the characteristic systems of any form
$\eta_\alpha$ which is of constant rank (see the proof of Theorem
5.5 in \cite{mag} for details). In the case of $3$-quasi-Sasakian
manifolds of rank $4l+1$ the characteristic systems of the forms
$\eta_\alpha$ define three integrable distributions, given by
$\langle \xi_\alpha,\xi_\beta\rangle\oplus{\mathcal E}^{4m}$, for
any $\alpha\ne\beta$.

Now, according to \cite{mag}, we define for each
$\alpha\in\left\{1,2,3\right\}$ two tensor fields of type $(1,1)$
$\psi_\alpha$ and $\theta_\alpha$ on $M$. We put, for a
$3$-quasi-Sasakian manifold of rank $4l+3$,
\begin{equation*}
 \psi_\alpha X=\left\{
             \begin{array}{ll}
               \phi_\alpha X, & \hbox{if $X\in \Gamma({\mathcal{E}}^{4l+3})$;}\\
               0, & \hbox{if $X\in \Gamma({\mathcal{E}}^{4m})$;}
             \end{array}
              \right.
\ \textrm{  } \
 \theta_\alpha X=\left\{
             \begin{array}{ll}
               0, & \hbox{if $X\in \Gamma({\mathcal{E}}^{4l+3})$;}\\
               \phi_\alpha X, & \hbox{if $X\in \Gamma({\mathcal{E}}^{4m})$,}\\
             \end{array}
              \right.
\end{equation*}
and for a $3$-quasi-Sasakian manifold of rank $4l+1$,
\begin{equation*}
 \psi_\alpha X=\left\{
             \begin{array}{ll}
               \phi_\alpha X, & \hbox{if $X\in \Gamma({\mathcal{E}}^{4l})$;}\\
               0, & \hbox{if $X\in \Gamma({\mathcal{E}}^{4m+3})$;}
             \end{array}
              \right.
\ \textrm{  } \
 \theta_\alpha X=\left\{
             \begin{array}{ll}
               0, & \hbox{if $X\in \Gamma({\mathcal{E}}^{4l})$;}\\
               \phi_\alpha X, & \hbox{if $X\in \Gamma({\mathcal{E}}^{4m+3})$.}\\
             \end{array}
              \right.
\end{equation*}
Note that, for each $\alpha\in\left\{1,2,3\right\}$ we have
$\phi_\alpha=\psi_\alpha+\theta_\alpha$.
We  have given two different definitions of $\psi_\alpha$ and
$\theta_\alpha$, depending on the two possible ranks (for each $l$)
that correspond to the two types of $3$-quasi-Sasakian manifolds. It
should be noted, however, that in both cases $\psi_\alpha$ and
$\theta_\alpha$ coincide on the horizontal subbundle $\mathcal H$.
Next, we define a new (pseudo-Riemannian, in general) metric
$\bar{g}$ on $M$ setting
\begin{equation*}
\bar{g}\left(X,Y\right)=\left\{
                          \begin{array}{ll}
                            -d\eta_\alpha\left(X,\phi_\alpha Y\right), & \hbox{for $X,Y\in\Gamma({\mathcal E}^{4l})$;} \\
                            g\left(X,Y\right), & \hbox{elsewhere.}
                          \end{array}
                        \right.
\end{equation*}
Note that this definition is well posed by virtue of (iv) of
Proposition \ref{lemmi}. The metric $\bar g$ is in fact a compatible
metric and $(\phi_\alpha,\xi_\alpha,\eta_\alpha,\bar g)$ is a normal
almost $3$-contact metric structure, in general
non-$3$-quasi-Sasakian (cf. \cite{mag}). Concerning the Levi Civita
connection $\bar\nabla$ of the metric $\bar g$ we prove the
following useful formula.

\begin{proposition}
With the notation above, one has in a $3$-quasi-Sasakian manifold
\begin{equation}\label{formulaproiettabilita}
\bar\nabla_{X}\xi_\alpha=-\psi_\alpha X
\end{equation}
for any $X\in\Gamma\left(\mathcal H\right)$ and
$\alpha\in\left\{1,2,3\right\}$.
\end{proposition}
\begin{proof}
In the case $c=0$ the result is an immediate consequence of Lemma
2.3 in \cite{tanno}. As for the case $c\neq 0$, using the same
Lemma, we have
\begin{equation}\label{formulablair}
\bar\nabla'_{\!\alpha}\xi_\alpha=-\psi'_\alpha,
\end{equation}
where
\begin{equation*}
 \psi'_\alpha X=\left\{
             \begin{array}{ll}
               \phi_\alpha X, & \hbox{if $X\in \Gamma({\mathcal{E}}^{4l}\oplus\left\langle\xi_\beta,\xi_\gamma\right\rangle)$;}\\
               0, & \hbox{if $X\in \Gamma({\mathcal{E}}^{4m}\oplus\langle\xi_\alpha\rangle)$}
             \end{array}
              \right.
\end{equation*}
and $\bar\nabla'_{\!\alpha}$ is the Levi Civita connection
associated to the  compatible metric $\bar g'_\alpha$ defined by
\begin{equation*}
\bar g'_\alpha\left(X,Y\right)=\left\{
                          \begin{array}{ll}
                            -d\eta_\alpha\left(X,\phi_\alpha Y\right), & \hbox{for $X,Y\in\Gamma({\mathcal E}^{4l}\oplus\left\langle\xi_\beta,\xi_\gamma\right\rangle)$;} \\
                            g\left(X,Y\right), & \hbox{elsewhere.}
                          \end{array}
                        \right.
\end{equation*}
Note that $\psi_\alpha=\psi'_\alpha$  on $\Gamma(\mathcal H)$.
Now,
 considering  $X\in\Gamma(\mathcal H)$, we prove that
\begin{equation}\label{uguaglianza}
\bar\nabla'_{\!\alpha}\vphantom{I}_X\xi_\alpha=\bar\nabla_{X}\xi_\alpha.
\end{equation}
It should be noted that the metric $\bar g'_\alpha$,
as well as $\bar g$, preserves the orthogonal decomposition
$TM={\mathcal H}\oplus{\mathcal V}$, whereas $\xi_1,\xi_2,\xi_3$
are only orthogonal and not orthonormal with respect to $\bar
g'_\alpha$: indeed $\bar
g'_\alpha(\xi_\beta,\xi_\beta)=\frac{c}{2}$. Then $\bar
g|_{{\mathcal H}\times{\mathcal H}}=\bar g'_\alpha|_{{\mathcal
H}\times{\mathcal H}}$ and $\bar g|_{{\mathcal H}\times{\mathcal
V}}=\bar g'_\alpha|_{{\mathcal H}\times{\mathcal V}}$. Now, in
order to prove \eqref{uguaglianza}, we show preliminarily that
$\bar\nabla'_{\!\alpha}\vphantom{I}_X\xi_\alpha, \bar\nabla_{X}\xi_\alpha \in
\Gamma\left(\mathcal H\right)$. Indeed,
\begin{align*}
2\bar g'_\alpha(\bar\nabla'_{\!\alpha}\vphantom{I}_X\xi_\alpha,\xi_\delta)&=X(\bar
g'_\alpha(\xi_\alpha,\xi_\delta))+\xi_\alpha(\bar
g'_\alpha(\xi_\delta,X))-\xi_\delta(\bar g'_\alpha(X,\xi_\alpha))\\
&\quad+\bar g'_\alpha([X,\xi_\alpha],\xi_\delta)+\bar
g'_\alpha([\xi_\delta,X],\xi_\alpha)-\bar
g'_\alpha([\xi_\alpha,\xi_\delta],X)=0,
\end{align*}
since $\bar g'_\alpha(\xi_\alpha,\xi_\delta)$ is constant,
$[\xi_\delta, \Gamma({\mathcal H})]\subset\Gamma({\mathcal H})$ for
any $\delta$, and $\mathcal V$ is integrable. Analogously, $\bar
g_\alpha(\bar\nabla_X\xi_\alpha,\xi_\delta)=0$. Then, using the
definitions of $\bar g$ and $\bar g'_\alpha$, we have that for any
$X,Y\in\Gamma\left(\mathcal H\right)$,
\begin{align*}
2\bar g(\bar\nabla_{X}\xi_\alpha-\bar\nabla'_{\!\alpha}\vphantom{I}_X\xi_\alpha,Y)
&=X(\bar g(\xi_\alpha,Y))+\xi_\alpha(\bar g(X,Y))-Y(\bar
g(\xi_\alpha,X))\\
&\quad+\bar g([X,\xi_\alpha],Y)+\bar g([Y,X],\xi_\alpha)-\bar g([\xi_\alpha,Y],X)\\
&\quad-X(\bar
g'_\alpha(\xi_\alpha,Y))-\xi_\alpha(\bar g'_\alpha(X,Y))+Y(\bar g'_\alpha(\xi_\alpha,X))\\
&\quad-\bar
g'_\alpha([X,\xi_\alpha],Y)-\bar g'_\alpha([Y,X],\xi_\alpha)+\bar
g'_\alpha([\xi_\alpha,Y],X)\\
&=\xi_\alpha(\bar g(X,Y))+\bar g([X,\xi_\alpha],Y)+ g([Y,X],\xi_\alpha)\\
&\quad-\bar g([\xi_\alpha,Y],X)-\xi_\alpha(\bar
g'_\alpha(X,Y))-\bar
g'_\alpha([X,\xi_\alpha],Y)\\
&\quad-g([Y,X],\xi_\alpha)+\bar g'_\alpha([\xi_\alpha,Y],X)=0.
\end{align*}
Therefore we have that $\bar\nabla_X \xi_\alpha=\bar\nabla'_{\!\alpha}\vphantom{I}_X\xi_\alpha
=-\psi'_\alpha X=-\psi_\alpha X$ and
\eqref{formulaproiettabilita} is proved.
\end{proof}

\begin{lemma}\label{lemmapro}
In any $3$-quasi-Sasakian manifold we have, for a cyclic permutation
$(\alpha,\beta,\gamma)$ of $\left\{1,2,3\right\}$,
\begin{equation}\label{fittizia}
{\mathcal L}_{\xi_\alpha}d\eta_\beta=cd\eta_\gamma.
\end{equation}
\end{lemma}
\begin{proof}
From the Cartan formula for the Lie derivative it follows that
${\mathcal
L}_{\xi_\alpha}d\eta_\beta=i_{\xi_\alpha}d^2\eta_\beta+di_{\xi_\alpha}d\eta_\beta=di_{\xi_\alpha}d\eta_\beta$,
so that it is enough to compute $i_{\xi_\alpha}d\eta_\beta$. By (i)
of Proposition \ref{lemmi} we have, for any
$X\in\Gamma\left(\mathcal H\right)$,
\[
(i_{\xi_\alpha}d\eta_\beta)(X)=2d\eta_\beta(\xi_\alpha,X)=0=c\eta_\gamma(X).\]
Now, distinguishing the cases $c=0$ and $c\ne0$, one can verify that
$i_{\xi_\alpha}d\eta_\beta=c\eta_\gamma$ also holds on $\Gamma({\mathcal V})$,
thus getting the result.
\end{proof}

\begin{lemma}\label{lemmapro3}
For any $X\in\Gamma\left(\mathcal H\right)$ and
$Y\in\Gamma({\mathcal E}^{4m})$ we have
$\left[X,Y\right]\in\Gamma({\mathcal H})$.
\end{lemma}
\begin{proof}
For any $\alpha\in\left\{1,2,3\right\}$ one has
$\eta_\alpha\left(\left[X,Y\right]\right)=-2d\eta_\alpha\left(X,Y\right)=\left(i_{Y}d\eta_\alpha\right)\left(X\right)=0$,
since $Y\in\Gamma({\mathcal E}^{4m})$. Hence
$\left[X,Y\right]\in\bigcap_{\alpha=1}^{3}\ker\left(\eta_\alpha\right)={\mathcal
H}$.
\end{proof}

\begin{lemma}\label{lemmapro4}
Let $(M^{4n+3},\phi_\alpha,\xi_\alpha,\eta_\alpha,g)$ be a
$3$-quasi-Sasakian manifold. Then the Reeb vector fields are
infinitesimal automorphisms with respect to the distributions
${\mathcal E}^{4l}$ and ${\mathcal E}^{4m}$.
\end{lemma}
\begin{proof}
Let us assume $c\neq 0$. Fixing $\alpha\in\left\{1,2,3\right\}$, by
Lemma 2.2 in \cite{tanno} we have that $[\xi_\alpha,\Gamma({\mathcal
E}^{4l})]\subset\Gamma( {\mathcal
E}^{4l}\oplus\langle\xi_\beta,\xi_\gamma\rangle)$ and
$[\xi_\alpha,\Gamma({\mathcal E}^{4m})]\subset\Gamma({\mathcal
E}^{4m})$. Then the result follows from (ii) of Proposition
\ref{lemmi}. Analogously one obtains the claim for $c=0$.
\end{proof}

\begin{proposition}\label{pro}
In any $3$-quasi-Sasakian manifold we have
\begin{equation}\label{proiezionephi}
{\mathcal L}_{\xi_\alpha}\phi_\beta=c\psi_\gamma,
\end{equation}
for any cyclic permutation $(\alpha,\beta,\gamma)$ of
$\left\{1,2,3\right\}$.
\end{proposition}
\begin{proof}
That \eqref{proiezionephi} holds on $\mathcal V$ follows immediately
from a direct computation and from the definitions of the tensors
$\psi_\alpha$. Next, for any $X\in\Gamma\left(\mathcal H\right)$ we
have
\begin{equation}\label{formulapro1}
\begin{split}
({\mathcal L}_{\xi_\alpha}\phi_\beta)X&=\left[\xi_\alpha,\phi_\beta
X\right]-\phi_\beta\left[\xi_\alpha,X\right]\\
&=\bar\nabla_{\xi_\alpha}\phi_\beta X-\bar\nabla_{\phi_\beta
X}\xi_\alpha-\phi_\beta\bar\nabla_{\xi_\alpha}X+\phi_\beta\bar\nabla_{X}\xi_\alpha\\
&=(\bar\nabla_{\xi_\alpha}\phi_\beta)X-\bar\nabla_{\phi_\beta
X}\xi_\alpha+\phi_\beta\bar\nabla_{X}\xi_\alpha.
\end{split}
\end{equation}
Now, using \eqref{olszakformula}, we compute
$\bar\nabla_{\xi_\alpha}\phi_\beta$. Taking into account the
normality of the structure $(\phi_\beta,\xi_\beta,\eta_\beta)$, (i)
of Proposition \ref{lemmi} and the horizontality of $X$ we have
\begin{equation}\label{richiamoutile1}
2\bar
g((\bar\nabla_{\xi_\alpha}\phi_\beta)X,Y)=3d\bar\Phi_\beta(\xi_\alpha,\phi_\beta
X,\phi_\beta Y)-3d\bar\Phi_\beta(\xi_\alpha,X,Y).
\end{equation}
If $Y=\xi_\delta$ for some $\delta\in\left\{1,2,3\right\}$,
by (ii) of Proposition \ref{lemmi} and the integrability of the
distribution spanned by $\xi_1, \xi_2, \xi_3$, we have
\begin{align*}
3d\bar\Phi_\beta(\xi_\alpha,X,\xi_\delta)&=\xi_\alpha(\bar\Phi_\beta(X,\xi_\delta))+X(\bar\Phi_\beta(\xi_\delta,\xi_\alpha))+\xi_\delta(\bar\Phi_\beta(\xi_\alpha,X))\\
&\quad-\bar\Phi_\beta([\xi_\alpha,X],\xi_\delta)-\bar\Phi_\beta([X,\xi_\delta],\xi_\alpha)-\bar\Phi_\beta([\xi_\delta,\xi_\alpha],X)=0,
\end{align*}
and, in the same way, we find
$3d\bar\Phi_\beta(\xi_\alpha,\phi_\beta X,\phi_\beta \xi_\delta)=0$,
so that $(\bar\nabla_{\xi_\alpha}\phi_\beta)X\in\Gamma(\mathcal H)$. Now we prove that
\begin{equation}\label{richiamoutile2}
\bar g((\bar\nabla_{\xi_\alpha}\phi_\beta)X,Y)=(c-2)\bar
g(\psi_\gamma X,Y)
\end{equation}
for every $X,Y\in\Gamma(\mathcal H)$. Indeed, by \eqref{richiamoutile1} we have
\begin{align}\label{equazionepro}
2\bar
g((\bar\nabla_{\xi_\alpha}\phi_\beta)X,Y)&=\xi_\alpha(\bar\Phi_\beta(\phi_\beta
X,\phi_\beta Y))+\phi_\beta X(\bar\Phi_\beta(\phi_\beta
Y,\xi_\alpha))\\
&\quad+\phi_\beta Y(\bar\Phi_\beta(\xi_\alpha,\phi_\beta \nonumber
X))-\bar\Phi_\beta([\xi_\alpha,\phi_\beta X],\phi_\beta Y)\\
&\quad-\bar\Phi_\beta([\phi_\beta X,\phi_\beta
Y],\xi_\alpha)-\bar\Phi_\beta([\phi_\beta Y,\xi_\alpha],\phi_\beta
\nonumber X)\\ \nonumber
&\quad-\xi_\alpha(\bar\Phi_\beta(X,Y))-X(\bar\Phi_\beta(Y,\xi_\alpha))-Y(\bar\Phi_\beta(\xi_\alpha,X))\\
\nonumber
&\quad+\bar\Phi_\beta([\xi_\alpha,X],Y)+\bar\Phi_\beta([X,Y],\xi_\alpha)+\bar\Phi_\beta([Y,\xi_\alpha],X)\\
\nonumber
&=\xi_\alpha(\bar\Phi_\beta(\phi_\beta X,\phi_\beta
\nonumber
Y))-\bar\Phi_\beta([\xi_\alpha,\phi_\beta X],\phi_\beta Y)\\
&\quad-\bar\Phi_\beta([\phi_\beta X,\phi_\beta \nonumber
Y],\xi_\alpha)-\bar\Phi_\beta([\phi_\beta Y,\xi_\alpha],\phi_\beta
X)\\ \nonumber
&\quad-\xi_\alpha(\bar\Phi_\beta(X,Y))+\bar\Phi_\beta([\xi_\alpha,X],Y)\\
&\quad+\bar\Phi_\beta([X,Y],\xi_\alpha)+\bar\Phi_\beta([Y,\xi_\alpha],X).
\nonumber
\end{align}
Because of the $\bar g$-orthogonal decomposition ${\mathcal
H}={\mathcal E}^{4l}\oplus{\mathcal E}^{4m}$ we can distinguish the
cases  $X,Y\in\Gamma({\mathcal E}^{4l})$, or $X\in\Gamma({\mathcal
E}^{4l})$, $Y\in\Gamma({\mathcal E}^{4m})$, or $X\in\Gamma({\mathcal
E}^{4m})$, $Y\in\Gamma({\mathcal E}^{4l})$, or
$X,Y\in\Gamma({\mathcal E}^{4m})$. In the first case, taking into
account that $\phi_\beta({\mathcal E}^{4l})={\mathcal E}^{4l}$ and
$[\xi_\alpha,\Gamma({\mathcal E}^{4l})]\subset\Gamma({\mathcal
E}^{4l})$ (cf. Lemma \ref{lemmapro4}) we get
\begin{align}\label{formulapro2}
2\bar
g((\bar\nabla_{\xi_\alpha}\phi_\beta)X,Y)&=({\mathcal L}_{\xi_\alpha}d\eta_\beta)(\phi_\beta X,\phi_\beta
Y)+\eta_\gamma([\phi_\beta X,\phi_\beta Y])\\ \nonumber
&\quad-({\mathcal
L}_{\xi_\alpha}d\eta_\beta)(X,Y)-\eta_\gamma([X,Y])\\
\nonumber &=({\mathcal L}_{\xi_\alpha}d\eta_\beta)(\phi_\beta
X,\phi_\beta
Y)-2d\eta_\gamma(\phi_\beta X,\phi_\beta Y)\\
&\quad-({\mathcal
L}_{\xi_\alpha}d\eta_\beta)(X,Y)+2d\eta_\gamma(X,Y). \nonumber
\end{align}
Continuing the computation and using  Lemma \ref{lemmapro} and (v) of
Proposition \ref{lemmi}, \eqref{formulapro2} becomes
\begin{align*}
2\bar g((\bar\nabla_{\xi_\alpha}\phi_\beta)X,Y)&=cd\eta_\gamma(\phi_\beta
X, \phi_\beta Y)-2d\eta_\gamma(\phi_\beta X, \phi_\beta
Y)-cd\eta_\gamma(X,Y)+2d\eta_\gamma(X,Y)\\
&=-c d\eta_\gamma(X,Y)+2d\eta_\gamma(X,Y)-c
d\eta_\gamma(X,Y)+2d\eta_\gamma(X,Y)\\
&=2(2-c)d\eta_\gamma(X,Y)\\
&=2(c-2)\bar g(\psi_\gamma X,Y).
\end{align*}
If we take $X\in\Gamma({\mathcal E}^{4l})$ and $Y\in\Gamma({\mathcal
E}^{4m})$ then, due to the orthogonality between ${\mathcal E}^{4l}$
and ${\mathcal E}^{4m}$, \eqref{equazionepro} reduces to
\begin{align*}
2\bar g((\bar\nabla_{\xi_\alpha}\phi_\beta)X,Y)&=-\bar\Phi_\beta([\phi_\beta
X,\phi_\beta Y],\xi_\alpha)+\bar\Phi_\beta([X,Y],\xi_\alpha)\\
&=\eta_\gamma([\phi_\beta X,\phi_\beta Y])-\eta_\gamma([X,Y])=0
\end{align*}
by Lemma \ref{lemmapro3}. Since $\bar g(\psi_\gamma X,Y)=0$, we get
\eqref{richiamoutile2}. Next, arguing as above, one finds that
\eqref{richiamoutile2} also holds  for $X\in\Gamma({\mathcal
E}^{4m})$ and $Y\in\Gamma({\mathcal E}^{4l})$.
 Finally, if $X,Y\in\Gamma({\mathcal E}^{4m})$, by the definition of $\bar g$ and  $d\Phi_\beta=0$, one has
\begin{align*}
2\bar
g((\bar\nabla_{\xi_\alpha}\phi_\beta)X,Y)&=3d\Phi_\beta(\xi_\alpha,\phi_\beta X,\phi_\beta
Y)-3d\Phi_\beta(\xi_\alpha,X,Y)=0
\end{align*}
which proves \eqref{richiamoutile2}, since $\psi_\gamma X=0$. Therefore we get  that
\begin{equation}\label{conclusione}
(\bar\nabla_{\xi_\alpha}\phi_\beta)X=(c-2)\psi_\gamma X
\end{equation}
for any $X\in\Gamma\left(\mathcal H\right)$. Continuing the
computation in \eqref{formulapro1},  we obtain, by virtue of
\eqref{conclusione} and \eqref{formulaproiettabilita},
\begin{equation*}
({\mathcal
L}_{\xi_\alpha}\phi_\beta)X=(c-2)\psi_\gamma
X+\psi_\alpha\phi_\beta X-\phi_\beta\psi_\alpha X,
\end{equation*}
so that $\left({\mathcal
L}_{\xi_\alpha}\phi_\beta\right)X=c\phi_\gamma X$ if
$X\in\Gamma({\mathcal E}^{4l})$ and $\left({\mathcal
L}_{\xi_\alpha}\phi_\beta\right)X=0$ if $X\in\Gamma({\mathcal
E}^{4m})$, from which the assertion follows.
\end{proof}

Now we prove the first strong consequence of Proposition \ref{pro},
namely that the only $3$-quasi-Sasakian manifolds whose Reeb vector
fields commute are the $3$-cosymplectic manifolds. We need the
following lemma.

\begin{lemma}
In any $3$-quasi-Sasakian manifold we have
\begin{equation}\label{formuladerivatelie}
({\mathcal L}_{\xi_\alpha}\Phi_\beta)(X,Y)=g(X,({\mathcal
L}_{\xi_\alpha}\phi_\beta)Y).
\end{equation}
for all $X,Y\in\Gamma\left(TM\right)$ and
$\alpha,\beta\in\{1,2,3\}$.
\end{lemma}
\begin{proof}
The assertion follows immediately from the fact that each
$\xi_\alpha$ is  Killing.
\end{proof}

\begin{theorem}\label{pro2}
Every $3$-quasi-Sasakian manifold of rank $4l+1$  is $3$-cosymplectic.
\end{theorem}
\begin{proof}
Using (iii) of Proposition \ref{lemmi} and
\eqref{formuladerivatelie}, we have
\begin{equation*}
2d\eta_\alpha(X,Y)=({\mathcal
L}_{\xi_\beta}\Phi_\gamma)(X,Y)=g(X,({\mathcal
L}_{\xi_\beta}\phi_\gamma)Y)
\end{equation*}
and the last term vanishes since, for $c=0$, ${\mathcal
L}_{\xi_\beta}\phi_\gamma=0$ due to Proposition \ref{pro}.
\end{proof}

\begin{corollary}\label{riccicurv1}
Any $3$-quasi-Sasakian manifold is Ricci-flat if and only if it is
$3$-cosymplectic.
\end{corollary}
\begin{proof}
That any $3$-cosymplectic manifold is Ricci-flat it has been
proven in \cite{cappellettidenicola}. Conversely, if a
$3$-quasi-Sasakian manifold
$(M,\phi_\alpha,\xi_\alpha,\eta_\alpha,g)$ is Ricci-flat, then by
\eqref{olszakricci} we get $\nabla\xi_\alpha=0$ for all
$\alpha\in\left\{1,2,3\right\}$, hence $c=0$. Thus applying
Theorem \ref{pro2} we get the result.
\end{proof}

It should be remarked that Corollary \ref{riccicurv1} may be
thought as an odd-dimensional analogue of the well-known fact that
any quaternionic-K\"{a}hler manifold is Ricci-flat if and only if
it is (locally) hyper-K\"{a}hler.

\section{$3$-quasi-Sasakian manifolds of rank $4l+3$}
We recall that an \emph{almost $\alpha$-Sasakian
manifold}\footnote{In the sequel by an abuse of notation we will use
the same symbol $\alpha$ both as an index and for indicating an
$\alpha$-Sasakian structure. The different meaning of this symbol
will be clear by the context.} (\cite{janssens}) is an almost
contact metric manifold satisfying $d\eta=\alpha\Phi$ for some
$\alpha\in\mathbb{R}^{\ast}$. An almost $\alpha$-Sasakian manifold
which is also normal is called an \emph{$\alpha$-Sasakian manifold}.
It is well known that an almost contact metric manifold is
$\alpha$-Sasakian if and only if
\begin{equation}\label{defalphasasaki}
(\nabla_{X}\phi)Y=\alpha(g(X,Y)\xi-\eta(Y)X)
\end{equation}
holds for all $X,Y\in\Gamma(TM)$, for some
$\alpha\in\mathbb{R}^{\ast}$. From \eqref{defalphasasaki} it
follows also that
\begin{equation}\label{conseguenze}
\nabla_{X}\xi=-\alpha\phi X, \ \
R(X,Y)\xi=\alpha^2(\eta(Y)X-\eta(X)Y).
\end{equation}
Since the fundamental $2$-form of an $\alpha$-Sasakian manifold is
exact (in particular closed) then the manifold is quasi-Sasakian.

Now consider an almost $3$-contact metric manifold
$(M,\phi_\delta,\xi_\delta,\eta_\delta,g)$ of dimension $4n+3$,
such that each structure is $\alpha$-Sasakian, and suppose $d\eta_\delta=\alpha_\delta \Phi_\delta$ for any
$\delta\in\left\{1,2,3\right\}$. Then, as it has been showed in
\cite{kashiwada2}, $\alpha_1=\alpha_2=\alpha_3:=a$, and we
have
\begin{equation}\label{relazionicampi}
[\xi_\alpha,\xi_\beta]=\nabla_{\xi_\alpha}\xi_\beta-\nabla_{\xi_\beta}\xi_\alpha=-a\phi_{\beta}\xi_\alpha+a\phi_\alpha\xi_\beta=2a\xi_\gamma.
\end{equation}
Hence, $M$ is $3$-quasi-Sasakian with $c=2a$ and maximal rank
$4n+3$. We will call an almost $3$-contact metric manifold such that
each structure is (almost) $\alpha$-Sasakian simply by a
\emph{(almost) $3$-$\alpha$-Sasakian manifold}.

An example of these manifolds is given by the sphere $S^{4n+3}(r)$
of radius $r$, considered as a hypersurface in $\mathbb{H}^{n+1}$.
Indeed, taking the quaternionic structure $(J_1,J_2,J_3)$ on
$\mathbb{H}^{n+1}$, one can define three vector fields on the
sphere, $\xi_\alpha=-J_\alpha \nu$, $\nu$ being a unit normal of
$S^{4n+3}(r)$. Next, one defines the tensor fields $\phi_\alpha$ of
type $(1,1)$ and the $1$-forms $\eta_\alpha$ by requiring that, for
any vector field $X$ tangent to the sphere, $\phi_\alpha X$ and
$\eta_\alpha(X)\nu$ are respectively the tangential and the normal
component of $J_\alpha X$ to the sphere. Considering the induced
Riemannian metric $g$, one obtains an almost $3$-contact metric
structure $(\phi_\alpha,\xi_\alpha,\eta_\alpha,g)$ which is
$3$-$\alpha$-Sasakian, since it is hyper-normal and the fundamental
$2$-forms satisfy $d\eta_\alpha=\frac{1}{r}\Phi_\alpha$.

We prove that, in fact, strictly almost $3$-$\alpha$-Sasakian
manifolds do not exist. This is a consequence of a generalization of
the Hitchin Lemma, due to Kashiwada, which we now recall.

\begin{lemma}[\cite{kashiwada1}] \label{genhitchin}
Let $(M^{4n},J_\alpha,G)$, $\alpha\in\left\{1,2,3\right\}$, be an
almost hyper-Hermitian manifold such that each fundamental
$2$-form $\Omega_\alpha$ satisfies
$d\Omega_\alpha=2\omega\wedge\Omega_\alpha$, for some $1$-form
$\omega$. Then each $J_\alpha$ is integrable.
\end{lemma}

\begin{proposition}\label{lemmamino}
Every almost $3$-$\alpha$-Sasakian manifold is necessarily
$3$-$\alpha$-Sasakian.
\end{proposition}
\begin{proof}
Let $(M,\phi_\alpha,\xi_\alpha,\eta_\alpha,g)$,
$\alpha\in\left\{1,2,3\right\}$, be an almost
$3$-$\alpha$-Sasakian manifold and let us consider on the product
manifold $M\times\mathbb R$ the almost Hermitian structures $J_1$,
$J_2$, $J_3$ defined by
\begin{equation*}
J_\alpha\left(X,f\frac{d}{dt}\right)=\left(\phi_\alpha
X-f\xi_\alpha,\eta_\alpha\left(X\right)\frac{d}{dt}\right),
\end{equation*}
for any vector field $X$ on $M$ and any smooth function $f$ on
$M\times\mathbb R$, where we have denoted by $t$ the global
coordinate on $\mathbb R$. A straightforward computation shows
that $J_\alpha J_\beta=-J_\beta J_\alpha=J_\gamma$ for an even
permutation $(\alpha,\beta,\gamma)$ of $\left\{1,2,3\right\}$.
Moreover, it is simple to check that the Riemannian metric
$G=g+dt^2$ is compatible with respect to the hyper-complex
structure $(J_1,J_2,J_3)$. Computing the expressions of the
fundamental $2$-forms we find
\begin{equation}\label{espressioni1}
\Omega_\alpha\left(X,Y\right)=\Phi_\alpha\left(X,Y\right), \ \
\Omega_\alpha\left(X,\frac{d}{dt}\right)=-\eta_\alpha\left(X\right),
\end{equation}
for all $X,Y\in\Gamma\left(TM\right)$ and
$\alpha\in\left\{1,2,3\right\}$. From \eqref{espressioni1} it
follows that
\begin{equation}\label{espressioni2}
d\Omega_\alpha\left(X,Y,Z\right)=d\Phi_\alpha\left(X,Y,Z\right)=0,
\ \
d\Omega_\alpha\left(X,Y,\frac{d}{dt}\right)=-\frac{2}{3}d\eta_\alpha\left(X,Y\right),
\end{equation}
for every $X,Y,Z\in\Gamma\left(TM\right)$. In particular, we have
that, for each $\delta\in\left\{1,2,3\right\}$,
$d\Omega_\delta=2\omega\wedge\Omega_\delta$, where $\omega=-\alpha
dt$. By Lemma \ref{genhitchin}, this concludes the proof.
\end{proof}

We will prove that every $3$-quasi-Sasakian manifold of maximal
rank is necessarily $3$-$\alpha$-Sasakian. This will be an
immediate consequence of the following result, which is an
analogue of Theorem \ref{pro2} for the class of $3$-quasi-Sasakian
manifolds which are not $3$-cosymplectic.

\begin{theorem}\label{pro3}
Let \ $(M^{4n+3},\phi_\alpha,\xi_\alpha,\eta_\alpha,g)$ \ be \ a \
$3$-quasi-Sasakian \ manifold \ such \ that \
$\left[\xi_\alpha,\xi_\beta\right]=c\xi_\gamma$, $c\neq 0$.
Let $4l+3$ be the rank of
$M^{4n+3}$. Then, for each $\alpha\in\left\{1,2,3\right\}$,
\begin{equation}\label{formulapesante1}
d\eta_\alpha(X,Y)=\frac{c}{2}g(X,\psi_\alpha Y)
\end{equation}
for all $X,Y\in\Gamma(TM)$. Consequently, on ${\mathcal
E}^{4l+3}$,
\begin{equation}\label{formulapesante}
d\eta_\alpha=\frac{c}{2}\Phi_\alpha.
\end{equation}
\end{theorem}
\begin{proof}
Using \eqref{formuladerivatelie} and \eqref{proiezionephi}, for
all $X,Y\in\Gamma(TM)$, we have
\begin{equation*}
({\mathcal L}_{\xi_\beta}\Phi_\gamma)(X,Y)=g(X,({\mathcal
L}_{\xi_\beta}\phi_\gamma)Y)=c g(X,\psi_\alpha Y).
\end{equation*}
On the other hand, by (iii) of Proposition \ref{lemmi}, ${\mathcal
L}_{\xi_\beta}\Phi_\gamma=2d\eta_\alpha$ from which
\eqref{formulapesante1} follows.
\end{proof}

\begin{corollary}\label{maximal}
\ Every \ $3$-quasi-Sasakian \ manifold \ of \ maximal \ rank \ is
\ necessarily \ $3$-$\alpha$-Sasakian.
\end{corollary}

\begin{remark}
It should be emphasized that in general no analogue of Theorem
\ref{pro3}, as well as of Corollary \ref{maximal},  holds for a
quasi-Sasakian manifold. These properties are thus a special feature
of $3$-quasi-Sasakian manifolds (cf. \cite{blair0,tanno}).
\end{remark}

\begin{corollary}
Let $(M,\phi_\alpha,\xi_\alpha,\eta_\alpha,g)$ be a
$3$-quasi-Sasakian manifold. Then for each
$\alpha\in\left\{1,2,3\right\}$
\begin{equation}\label{formulapesante2}
\nabla\xi_\alpha=-\frac{c}{2}\psi_\alpha.
\end{equation}
\end{corollary}
\begin{proof}
By Theorem \ref{pro2} the assertion easily holds for $c=0$, so
that we can assume $c\neq 0$. Then by \eqref{olszakformula} we
have
$g((\nabla_X\phi_\alpha)\xi_\alpha,Z)=d\eta_\alpha(\phi_\alpha
\xi_\alpha,X)\eta_\alpha(Z)-d\eta_\alpha(\phi_\alpha
Z,X)\eta_\alpha(\xi_\alpha)=-d\eta_\alpha(\phi_\alpha Z,X)$, from
which, applying \eqref{formulapesante1}, it follows that
\begin{equation*}
g(\nabla_{X}\xi_\alpha,\phi_\alpha Z)=-d\eta_\alpha(\phi_\alpha
Z,X)=-\frac{c}{2}g(\psi_\alpha X,\phi_\alpha Z).
\end{equation*}
Therefore, tacking into account the fact that
$g(\nabla_{X}\xi_\alpha,\xi_\alpha)=0=-\frac{c}{2}g(\psi_\alpha
X,\xi_\alpha)$, \eqref{formulapesante2} is proved.
\end{proof}

\begin{corollary}\label{corollariopesante}
In any $3$-quasi-Sasakian manifold one has, for each
$\alpha\in\left\{1,2,3\right\}$,
\begin{equation}\label{altraformulapesante}
\left(\nabla_{X}\phi_\alpha\right)Y=\frac{c}{2}\left(\eta_\alpha\left(Y\right)\psi_\alpha^2
X -g\left(\psi_\alpha^2
X,
Y\right)\xi_\alpha\right),
\end{equation}
for any $X,Y\in\Gamma(TM)$. Consequently, for any $X\in\Gamma(TM)$ and $\alpha\ne\beta$,
\begin{equation}\label{nablaphi}
(\nabla_{\xi_\beta}\phi_\alpha)X=\frac{c}{2}(\eta_\beta(X)\xi_\alpha-\eta_\alpha(X)\xi_\beta).
\end{equation}
\end{corollary}
\begin{proof}
It is a consequence of \eqref{formulaquasisasaki}, \eqref{formulapesante2} and the fact that $\phi_\alpha\psi_\alpha=\psi_\alpha^2$.
\end{proof}

\begin{corollary}\label{curvatura4l+3}
In any $3$-quasi-Sasakian manifold of rank $4l+3$ the Riemannian curvature satisfies:
\begin{enumerate}
\item[(i)] $R(X,Y)\xi_\alpha=0$ for any $X,Y\in \Gamma({\mathcal E}^{4m})$ and $\alpha\in\{1,2,3\}$;
\item[(ii)] $R(X,\xi_\beta)\xi_\alpha=0$ for any $X\in \Gamma({\mathcal E}^{4m})$ and $\alpha,\beta\in\{1,2,3\}$;
\item[(iii)] $R(X,\xi_\alpha)\xi_\alpha=\frac{c^2}{4}X$ for any $X\in \Gamma({\mathcal E}^{4l})$ and $\alpha\in\{1,2,3\}$;
\item[(iv)] $R(X,\xi_\beta)\xi_\alpha=0$ for any $X\in \Gamma({\mathcal E}^{4l})$ and $\alpha\ne\beta$.
\end{enumerate}
Consequently, for the sectional curvatures we have $K(\xi_\alpha,X)=0$ for any $X\in\Gamma({\mathcal E}^{4m})$ and $K(\xi_\alpha,X)=\frac{c^2}{4}\|X\|^2$ for any $X\in\Gamma({\mathcal E}^{4l})$.
\end{corollary}
\begin{proof}
All formulas are obtained by direct computations using \eqref{formulapesante2} and the definition of $\psi_\alpha$. In (i) we also apply the integrability of the distribution ${\mathcal E}^{4m}$, while in (ii) we use the fact that each $\xi_\alpha$ is an
infinitesimal automorphism with respect to the distributions
${\mathcal E}^{4m}$. As for (iv), a direct computation shows that $R(X,\xi_\beta)\xi_\alpha=\frac{c}{2}(\nabla_{\xi_\beta}\phi_\alpha)X$ for any $X\in \Gamma({\mathcal E}^{4l})$ and $\alpha\ne\beta$. Applying \eqref{nablaphi}, we get the result.
\end{proof}

\section{Transverse geometry of a $3$-quasi-Sasakian manifold}

In this section we prove that every $3$-quasi-Sasakian manifold is
canonically foliated by four Riemannian and totally geodesic
foliations, which play a fundamental role in the study of
$3$-quasi-Sasakian geometry. In particular we study the leaf spaces
of these foliations which, as we will see, greatly influence the
Riemannian geometry of the $3$-quasi-Sasakian manifold in question.


We start with the study of the $3$-dimensional foliation $\mathcal
V$ defined by the Reeb vector fields. Let $\nabla^B$ be the Bott
connection associated to $\mathcal V$, that is the partial
connection on the normal bundle $TM/{\mathcal V}\cong\mathcal H$ of
$\mathcal V$ defined by
\begin{equation*}
\nabla^B_{V}Z:=[V,Z]_{\mathcal H}
\end{equation*}
for all $V\in\Gamma(\mathcal V)$ and $Z\in\Gamma(\mathcal H)$,
where the subscript $\mathcal H$ denotes  the horizontal component.
Following \cite{tondeur} we may construct an adapted connection on
$\mathcal H$ putting
\begin{equation*}
\tilde\nabla_{X}Y:=\left\{
                     \begin{array}{ll}
                       \nabla^B_{X}Y, & \hbox{if $X\in\Gamma(\mathcal V)$;} \\
                       (\nabla_{X}Y)_{\mathcal H}, & \hbox{if $X\in\Gamma(\mathcal H)$.}
                     \end{array}
                   \right.
\end{equation*}
This connection can be also extended to a connection on all $TM$ by
requiring that $\tilde\nabla\xi_\alpha=0$ for each
$\alpha\in\left\{1,2,3\right\}$. Some properties of this global
connection have been considered in \cite{cappellettidenicola} for
any almost $3$-contact metric manifold. Now combining Theorem
\ref{principale} with Theorem 3.6 in \cite{cappellettidenicola} we
have:

\begin{theorem}\label{connessionecanonica}
Let $(M,\phi_\alpha,\xi_\alpha,\eta_\alpha,g)$ be a
$3$-quasi-Sasakian manifold. Then there exists a unique connection
$\tilde\nabla$ on $M$ satisfying the following properties:
\begin{enumerate}
  \item[(i)]
  $\tilde\nabla\eta_\alpha=0$, $\tilde\nabla\xi_\alpha=0$, for each
  $\alpha\in\left\{1,2,3\right\}$,
  \item[(ii)] $\tilde\nabla g=0$,
  \item[(iii)] $\tilde
  T\left(X,Y\right)=2\sum_{\alpha=1}^{3}d\eta_{\alpha}(X,Y)\xi_\alpha$,
  for all $X,Y\in\Gamma\left(TM\right)$.
\end{enumerate}
Furthermore, we have, for any cyclic permutation $(\alpha, \beta,
\gamma)$ of $\left\{1,2,3\right\}$,
\begin{equation*}
(\tilde\nabla_{X}\phi_\alpha)Y=
-c\left(\eta_\beta(X)\psi_\gamma(Y_{\mathcal
H})-\eta_\gamma(X)\psi_\beta(Y_{\mathcal H})\right).
\end{equation*}
\end{theorem}
\begin{proof}
Theorem 3.6 of \cite{cappellettidenicola} guarantees the existence
and the uniqueness of a linear connection $\tilde\nabla$ on $M$ such
that $\tilde\nabla\xi_\alpha=0$, $(\tilde\nabla_{Z}g)(X,Y)=0$ for
all $X,Y,Z\in\Gamma(\mathcal H)$ and $\tilde
T(X,Y)=2\sum_{\alpha=1}^{3}d\eta_\alpha(X,Y)\xi_\alpha$, $\tilde
T(X,\xi_\alpha)=0$ for all $X,Y\in\Gamma(\mathcal H)$. This
connection is explicitly defined as above. Since each $\xi_\alpha$
is Killing we have that $\tilde\nabla$ is metric
(\cite{cappellettidenicola}). Moreover, (i) of Proposition
\ref{lemmi} implies that each $1$-form $\eta_\alpha$ is
$\tilde\nabla$-parallel and, for the torsion tensor field, $\tilde
T(X,\xi_\alpha)=0=\sum_{\delta=1}^{3}2d\eta_\delta(X,\xi_\alpha)\xi_\delta$
for any $X\in\Gamma({\mathcal H})$ (cf. \cite{cappellettidenicola}).
Finally, from the integrability of $\mathcal V$ it follows also
$\tilde
T(\xi_\alpha,\xi_\beta)=[\xi_\beta,\xi_\alpha]=\sum_{\delta=1}^{3}2d\eta_\delta(\xi_\alpha,\xi_\beta)\xi_\delta$.
It remains to check the final part of the statement. We prove that
\begin{equation}\label{dimostrato}
(\tilde\nabla_{X}\phi_\alpha)Y=\left\{
\begin{array}{ll}
  0, & \hbox{for $X\in\Gamma\left({\mathcal H}\right)$ or $X=\xi_\alpha$ or $Y\in \Gamma\left({\mathcal V}\right)$;}\\
  -c\psi_\gamma Y, & \hbox{for  $X=\xi_\beta$, $Y\in\Gamma\left(\mathcal H\right)$ and $(\alpha, \beta, \gamma)$ cyclic.}
\end{array}
\right.
\end{equation}
Firstly, since the Reeb vector fields  are parallel with respect to
$\tilde{\nabla}$,  one has $(\tilde{\nabla}_{X}\phi_1)\xi_\alpha$
$=0$ for any $\alpha\in\left\{1,2,3\right\}$. Next, taking
$X,Y\in\Gamma\left(\mathcal H\right)$, using
\eqref{formulaquasisasaki} we have
\begin{align*}
(\tilde\nabla_{X}\phi_1)Y&=(\nabla_{X}\phi_1
Y)_{\mathcal H}-\phi_1(\nabla_{X}Y)_{\mathcal H}\\
&=\nabla_{X}\phi_1 Y-\sum_{\alpha=1}^{3}g(\nabla_{X}\phi_1
Y,\xi_\alpha)\xi_\alpha-\phi_1\nabla_{X}Y+\sum_{\alpha=1}^{3}g(\nabla_{X}Y,\xi_\alpha)\phi_1\xi_\alpha\\
&=(\nabla_{X}\phi_1)Y+\sum_{\alpha=1}^{3}g(\phi_1
Y,\nabla_{X}\xi_\alpha)\xi_\alpha+g(\nabla_{X}Y,\xi_2)\xi_3-g(\nabla_{X}Y,\xi_3)\xi_2\\
&=-g(\nabla_{X}\xi_1,\phi_1
Y)\xi_1-\eta_1(Y)\phi_1\nabla_{X}\xi_1+g(\phi_1
Y,\nabla_{X}\xi_1)\xi_1\\
&\quad+g(\phi_1 Y,\nabla_{X}\xi_2)\xi_2+g(\phi_1
Y,\nabla_{X}\xi_3)\xi_3-g(Y,\nabla_{X}\xi_2)\xi_3+g(Y,\nabla_{X}\xi_3)\xi_2\\
&=-g(\phi_1\nabla_{X}\xi_3+\nabla_{X}\xi_2,Y)\xi_3-g(\phi_1\nabla_{X}\xi_2-\nabla_{X}\xi_3,Y)\xi_2=0.
\end{align*}
Indeed, one has \(
\nabla_{X}\xi_2=-\nabla_{X}\phi_1\xi_3=-\phi_1\nabla_{X}\xi_3, \)
since \eqref{formulaquasisasaki} and the facts that $\xi_1$ is
Killing and $\mathcal V$ is totally geodesic imply
$(\nabla_{X}\phi_1)\xi_3=g(\nabla_{X}\xi_1,\xi_2)\xi_1 = -g(\nabla_{\xi_2}\xi_1,X)\xi_1$ $=0$.
Analogously, $\nabla_{X}\xi_3=\phi_1\nabla_{X}\xi_2$. Finally, for
any $Y\in\Gamma\left(\mathcal H\right)$, by the definition of
$\tilde{\nabla}$ and by (ii) of Proposition \ref{lemmi} one has
\begin{equation*}
(\tilde{\nabla}_{\xi_1}\phi_1)Y=\nabla^B_{\xi_1}\phi_1
Y-\phi_1\nabla^B_{\xi_1}Y=\left[\xi_1,\phi_1
Y\right]-\phi_1\left[\xi_1,Y\right]=\left({\mathcal
L}_{\xi_1}\phi_1\right)Y=0.
\end{equation*}
Similarly, using also Proposition \ref{pro} we have
$(\tilde{\nabla}_{\xi_2}\phi_1)Y=\left({\mathcal
L}_{\xi_2}\phi_1\right)Y=-c\psi_3 Y$ and
$(\tilde{\nabla}_{\xi_3}\phi_1)Y=\left({\mathcal
L}_{\xi_3}\phi_1\right)Y=c\psi_2 Y$. We have thus proved
\eqref{dimostrato}. Now, decomposing any vector fields $X$ and $Y$
in their horizontal and vertical components one easily gets the
claimed formula for $\tilde{\nabla}\phi_\alpha$.
\end{proof}

%
%
%
%

Using the constructions above, we prove the projectability of a
$3$-quasi-Sasakian structure. Indeed we know by Theorem
\ref{principale} that a $3$-quasi-Sasakian manifold $M$ of dimension
$4n+3$ is foliated by a $3$-dimensional foliation $\mathcal V$
which, as we have seen, influences greatly the geometry of $M$. It
can be very useful to know more about the space of leaves
$M'=M/{\mathcal V}$ generated by this foliation, which is, under
some assumptions of regularity, a $4n$-dimensional smooth manifold,
more in general an orbifold. As $\mathcal V$ is a Riemannian
foliation, the metric $g$ projects along the leaves onto a
Riemannian metric $g'$ on $M'$. What we have to study is the (local)
projectability of the tensor fields $\phi_\alpha$ or, more in
general, of the subbundle of $\textrm{End}(TM)$ that they span. This
question is solved by the following Theorem.

\begin{theorem}\label{pro1}
Every $3$-quasi-Sasakian manifold admits a canonical, projectable,
transversal almost quaternionic-Hermitian structure.
\end{theorem}
\begin{proof}
Let $(M^{4n+3},\phi_\alpha,\xi_\alpha,\eta_\alpha,g)$ be a
$3$-quasi-Sasakian manifold of rank $4l+3$. We first notice that the
distributions ${\mathcal E}^{4l}$ and ${\mathcal E}^{4m}$ are
foliated objects, i.e. they locally project along the leaves of
$\mathcal V$ onto two distributions on the space of leaves which we
denote by ${\mathcal E'}^{4l}$ and ${\mathcal E'}^{4m}$,
respectively. In order to prove this, let $\pi$ be a local
submersion defining the foliation $\mathcal V$. We note that, as
$\mathcal V$ is a Riemannian foliation, $\pi$ is in fact a
Riemannian submersion. We have to prove that, given any two points
$x$ and $y$ on the same fiber, one has
\begin{equation}\label{benposta}
\pi_{\ast_x}({\mathcal E}^{4l}_{x})=\pi_{\ast_y}({\mathcal
E}^{4l}_{y}), \ \ \pi_{\ast_x}({\mathcal
E}^{4m}_{x})=\pi_{\ast_y}({\mathcal E}^{4m}_{y}).
\end{equation}
Firstly observe that from Lemma \ref{lemmapro4} it follows
immediately that the Bott connection preserves the distributions
${\mathcal E}^{4l}$ and ${\mathcal E}^{4m}$. In particular these
distributions are preserved by the parallel transport along vertical
curves. Now, let $x, y \in M$ such that $\pi(x)=x'=\pi(y)$ and let
$\gamma$ be a leaf curve such that $\gamma(0)=x$ and $\gamma(1)=y$.
Let $\tau$ denote the parallel transport with respect to the Bott
connection along the curve $\gamma$. Then we preliminarily prove
that the following diagram commutes,
$$
\xymatrix{{\mathcal E}_x^{4l}  \ar[rr]^{\tau} \ar[dr]_{\pi_{\ast_x}}&& {\mathcal E}_y^{4l} \ar[dl]^{\pi_{\ast_y}}\\
                    & T_{x'}M'&}
$$
Indeed, let $v\in{\mathcal E}_x^{4l}\subset{\mathcal H}_x$ and $X:I\rightarrow \mathcal
H$ be the unique vector field along $\gamma$ such that
$\nabla^B_{\gamma'}X\equiv 0$ and $X(0)=v$, so that $\tau(v)=X(1)$.
Let $Y'$ be any vector field on the base space and $Y$ be the
corresponding basic vector field on $M$. Then we have
\begin{align*}
\frac{d}{dt}g'(\pi_{\ast_{\gamma(t)}}(X(t)),Y'_{\pi(\gamma(t))})&=\frac{d}{dt}g'(\pi_{\ast_{\gamma(t)}}(X(t)),\pi_{\ast_{\gamma(t)}}(Y_{\gamma(t)}))\\
&=\frac{d}{dt}g(X(t),Y_{\gamma(t)})\\
&=g(\tilde\nabla_{\gamma'}X,Y)_{\gamma(t)}+g(X,\tilde\nabla_{\gamma'}Y)_{\gamma(t)}\\
&=g(\nabla^B_{\gamma'}X,Y)_{\gamma(t)}+g(X,\nabla^B_{\gamma'}Y)_{\gamma(t)}\\
&=g(X,\nabla^B_{\gamma'}Y)_{\gamma(t)}.
\end{align*}
Now, as  $\gamma'(t)\in{\mathcal V}_{\gamma(t)}$ for all $t\in I$,
$\gamma'=\sum_{\alpha=1}^{3}f_\alpha\xi_\alpha$  for some functions
$f_\alpha$, hence
\begin{equation*}
\nabla^B_{\gamma'}Y=\sum_{\alpha=1}^{3}f_\alpha\nabla^B_{\xi_\alpha}Y=\sum_{\alpha=1}^{3}f_\alpha[\xi_\alpha,Y]_{\mathcal
H}=0
\end{equation*}
because $Y$ is assumed to be basic. Therefore
\begin{equation*}
g'_{\pi(\gamma(0))}(\pi_{\ast_{\gamma(0)}}(X(0)),Y'_{\pi(\gamma(0))})=g'_{\pi(\gamma(1))}(\pi_{\ast_{\gamma(1)}}(X(1)),Y'_{\pi(\gamma(1))}),
\end{equation*}
that is
\begin{equation*}
g'_{x'}(\pi_{\ast_{x}}(v),Y'_{x})=g'_{x'}(\pi_{\ast_{y}}(\tau(v)),Y'_{x'}).
\end{equation*}
By the arbitrariness of $Y'$ we conclude that
$\pi_{\ast_x}(v)=(\pi_{\ast_y}\circ\tau)(v)$. Thus
$\pi_{\ast_x}({\mathcal E}^{4l}_{x})=\pi_{\ast_y}(\tau({\mathcal
E}^{4l}_{x}))=\pi_{\ast_y}({\mathcal E}^{4l}_{y})$ and arguing
analogously for ${\mathcal E}^{4m}$ one has $\pi_{\ast_x}({\mathcal
E}^{4m}_{x})=\pi_{\ast_y}({\mathcal E}^{4m}_{y})$. Hence
\eqref{benposta} are proved and ${\mathcal E}^{4l}$, ${\mathcal
E}^{4m}$ project to well-defined distributions ${\mathcal E}'^{4l}$,
${\mathcal E}'^{4m}$ which are also mutually orthogonal since the
Riemannian metric $g$ is bundle-like. We can now construct an almost
quaternionic structure on the space of leaves $M'$. By an abuse of
notation we will denote, for each $\alpha\in\left\{1,2,3\right\}$,
by $\psi_\alpha$ and $\theta_\alpha$ the restriction of
$\phi_\alpha$ to ${\mathcal E}^{4l}$ and ${\mathcal E}^{4m}$,
respectively. Let ${\bar{Q}}$ be the subbundle of
$\textrm{End}({\mathcal E}^{4l})$ spanned by $\psi_1, \psi_2,
\psi_3$ and $\bar{\bar{Q}}$ be the subbundle of
$\textrm{End}({\mathcal E}^{4m})$ spanned by $\theta_1, \theta_2,
\theta_3$. For any $X\in{\mathcal E}^{4l}$ we have
\begin{equation*}
(\nabla^B_{\xi_\alpha}\psi_\beta)X=[\xi_\alpha,\psi_\beta
X]_{\mathcal H}-\psi_\beta[\xi_\alpha,X]_{\mathcal
H}=({\mathcal
L}_{\xi_\alpha}\phi_\beta)X=c\psi_\gamma X,
\end{equation*}
by Lemma \ref{lemmapro4} and Proposition \ref{pro}. Thus the Bott
connection preserves the subbundle $\bar Q$ and this guarantees the
projectability of $\bar Q$  onto an almost quaternionic structure
$\bar{Q}'\subset \textrm{End}({\mathcal E}'^{4l})$ on the space of
leaves of the foliation $\mathcal V$ (cf. \cite{piccinnivaisman}).
For the subbundle $\bar{\bar Q}$ we can prove something more, namely
that each $\theta_\alpha$ is projectable. Indeed, for any
$Y\in{\mathcal E}^{4m}$ we have
\begin{equation*}
({\mathcal L}_{\xi_\alpha}\theta_\beta)Y=[\xi_\alpha,\phi_\beta
Y]-\phi_\beta[\xi_\alpha,Y]=({\mathcal
L}_{\xi_\alpha}\phi_\beta)Y=c\psi_\gamma Y=0
\end{equation*}
again by Lemma \ref{lemmapro4} and Proposition \ref{pro}. Thus each
$\theta_\alpha$ projects to a tensor field $\theta'_\alpha$ defined
on ${\mathcal E}'^{4m}$. Let us denote by $\bar{\bar{Q}}'$ the
subbundle of $\textrm{End}({\mathcal E}'^{4m})$ that is spanned by
$\theta'_1, \theta'_2, \theta'_3$. Since $TM'={\mathcal
E}'^{4l}\oplus{\mathcal E}'^{4m}$, from $\bar Q'$ and $\bar{\bar
Q}'$ we can define an almost quaternionic structure on $M'$ in the
following way. Let $\psi'_1, \psi'_2, \psi'_3$ be a local basis for
$\bar Q'$ defined on an open coordinate neighborhood $U'$. Then we
define three tensor fields, on $U'$, by
\begin{equation*}
\phi'_\alpha :=\left\{
                   \begin{array}{ll}
                     \psi'_\alpha, & \hbox{on ${\mathcal E}'^{4l}$,} \\
                     \theta'_\alpha, & \hbox{on ${\mathcal E}'^{4m}$,}
                   \end{array}
                 \right.
\end{equation*}
for each $\alpha\in\left\{1,2,3\right\}$. Let $Q'$ be the subbundle
of $\textrm{End}(TM')$ spanned by $\phi'_1, \phi'_2, \phi'_3$. Since
 in the overlapping of two coordinate neighborhoods $U'$ and
$V'$ the matrix of the components of the $\phi'_\alpha|_{U'}$ with
respect to the $\phi'_\alpha|_{V'}$ has the form
\begin{equation*}
\left(
  \begin{array}{cc}
    A & 0 \\
    0 & I_{4m} \\
  \end{array}
\right)
\end{equation*}
for some $A\in \emph{SO}(4l)$, we conclude that $Q'$ defines an
almost quaternionic-Hermitian structure on $M'$.
\end{proof}

We now examine more explicitly the case when a $3$-quasi-Sasakian
structure projects (locally) onto a hyper-K\"{a}hler structure.

\begin{proposition}\label{proiettabilita0}
Let $(M,\phi_\alpha,\xi_\alpha,\eta_\alpha,g)$ be a
$3$-quasi-Sasakian manifold of rank $4l+3$. Then the following statements are equivalent:
\begin{enumerate}
  \item[(i)] each structure tensor $\phi_\alpha$ is projectable with respect to  $\mathcal V$;
  \item[(ii)] for all $\alpha\in\left\{1,2,3\right\}$ $d\eta_\alpha=0$ on $\mathcal H$;
  \item[(iii)] the horizontal subbundle $\mathcal H$ is integrable;
  \item[(iv)] $M$ is locally a Riemannian product of a hyper-K\"{a}hler
  manifold and the $3$-dimensional sphere of constant curvature $\frac{c^2}{4}$;
  \item[(v)] $M$ is a $3$-quasi-Sasakian manifold of rank $3$.
\end{enumerate}
Furthermore, if  one of the above conditions holds, then the Ricci tensor of
$M$ is given by
\begin{equation}\label{primorisultatoricci}
\emph{Ric}=\frac{c^2}{2}\left(\eta_1\otimes\eta_1+\eta_2\otimes\eta_2+\eta_3\otimes\eta_3\right),
\end{equation}
hence $M$ is $\eta$-Einstein.
\end{proposition}
\begin{proof}
Each $\phi_\alpha$ is projectable if and only if  $({\mathcal
L}_{\xi_\beta}\phi_\alpha)X=0$ for all
$\beta\in\left\{1,2,3\right\}$ and $X\in\Gamma\left(\mathcal
H\right)$, which, by virtue of \eqref{formuladerivatelie}, is
equivalent to the vanishing of $({\mathcal
L}_{\xi_\beta}\Phi_\alpha)(X,Y)$ for all
$\beta\in\left\{1,2,3\right\}$ and $X,Y\in\Gamma\left(\mathcal
H\right)$.  This proves that (i) is equivalent to (ii). The
equivalence of (ii), (iii) and (v) is obvious. In order to prove the
equivalence of (iii) and (iv), let us assume the integrability of
$\mathcal H$. Then, since $\mathcal V$ defines a Riemannian
foliation with totally geodesic leaves, $M$ is a local Riemannian
product of a manifold $M'^{4n}$ tangent to the distribution
$\mathcal H$ and a $3$-dimensional manifold tangent to $\mathcal V$,
which is a space of constant curvature $\frac{c^2}{4}$ by virtue of
\eqref{Rxi}. We denote by $G$ the projection of the Riemannian
metric $g$ on $M'^{4n}$. The tensor fields $\phi_1$, $\phi_2$,
$\phi_3$ project to three tensor fields $J_1$, $J_2$, $J_3$ on
$M'^{4n}$ and it is easy to check that they satisfy the quaternionic
relations. In fact $\left(J_\alpha, G\right)$ are almost K\"{a}hler
structures which are integrable because $N^{(1)}_\alpha=0$. Now,
since $M$ is a local Riemannian product of a hyper-K\"{a}hler
manifold, which is Ricci-flat, and the sphere
$S^3\big(\frac{c^2}{4}\big)$ which is parallelized by the vector
fields $\xi_1$, $\xi_2$, $\xi_3$, the Ricci tensor of $M$ is given
by \eqref{primorisultatoricci}, and $M$ is $\eta$-Einstein.
\end{proof}

Now we study the transverse geometry of a $3$-quasi-Sasakian
manifold with respect to three other distinct foliations: ${\mathcal
E}^{4m}$, ${\mathcal E}^{4m+3}$ and ${\mathcal E}^{4l}$. We will
prove that ${\mathcal E}^{4m}$ is transversely
$3$-$\alpha$-Sasakian, ${\mathcal E}^{4m+3}$ transversely
quaternionic-K\"{a}hler and ${\mathcal E}^{4l}$ transversely
hyper-K\"{a}hler.

\begin{theorem}\label{integrabilita}
Let $(M,\phi_\alpha,\xi_\alpha,\eta_\alpha,g)$ be a
$3$-quasi-Sasakian manifold of rank $4l+3$. Then the integrable
distribution ${\mathcal E}^{4m}$ defines a Riemannian foliation of
dimension $4m$ whose leaves are hyper-K\"{a}hler manifolds.
Furthermore, the space of leaves of this foliation is
$3$-$\alpha$-Sasakian.
\end{theorem}
\begin{proof}
Let $N$ be a leaf of the foliation defined by ${\mathcal E}^{4m}$
and let, for each $\alpha\in\left\{1,2,3\right\}$, $J_\alpha$,
$\Omega_\alpha$, $G$ be the tensors on $N$ obtained from
$\phi_\alpha$, $\Phi_\alpha$, $g$ by restriction. Then
$(J_\alpha,\Omega_\alpha,G)$ defines an almost hyper-Hermitian
structure on $N$ which is integrable because  its Nijenhuis tensor
satisfies
$[J_\alpha,J_\alpha]=[\phi_\alpha,\phi_\alpha]|_N=([\phi_\alpha,\phi_\alpha]+2d\eta_\alpha\otimes\xi_\alpha)|_N=0$,
since $M$ is hyper-normal. We prove that the foliation ${\mathcal
E}^{4m}$ is transversely $3$-$\alpha$-Sasakian. We begin observing
that, for each $\alpha\in\left\{1,2,3\right\}$, the forms
$\eta_\alpha$ and $d\eta_\alpha$ are projectable, since for all
$V\in\Gamma({\mathcal E}^{4m})$ we have $i_{V}\eta_\alpha=0$ and
$i_{V}d\eta_\alpha=0$ by definition of ${\mathcal E}^{4m}$. Next, by
Lemma \ref{lemmapro4}, the Reeb vector fields $\xi_1$, $\xi_2$,
$\xi_3$ are basic vector fields. More delicate is the projectability
of the tensor fields $\phi_\alpha$. First note that as each $2$-form
$d\eta_\alpha$ is non-degenerate on ${\mathcal E}^{4l}$, it induces
a musical isomorphism $(d\eta_\alpha)^{\flat}:X \mapsto
d\eta_\alpha(X,\cdot)$ between ${\mathcal E}^{4l}$ and the $1$-forms
which vanish on ${\mathcal E}^{4m+3}$. We denote its inverse by
$(d\eta_\alpha)^{\sharp}$. As in Lemma 4.1 in
\cite{cappellettidenicola}, we have
\begin{equation}\label{formulaimportante}
\phi_\alpha X = -(d\eta_\beta)^{\sharp}(d\eta_\gamma)^{\flat}(X)
\end{equation}
for all $X\in\Gamma({\mathcal E}^{4l})$, where as usual
$(\alpha,\beta,\gamma)$ is an even permutation of
$\left\{1,2,3\right\}$. In order to prove that each $\phi_\alpha$ is
foliate with respect to the foliation ${\mathcal E}^{4m}$, it is
sufficient to show that $\phi_\alpha$ maps basic vector fields to
basic vector fields. In view of \eqref{formulaimportante} we prove
in fact that for each $\delta\in\left\{1,2,3\right\}$
$(d\eta_\delta)^{\flat}$ (respectively, $(d\eta_\delta)^{\sharp}$)
maps basic vector fields (respectively, basic $1$-forms) to basic
$1$-forms (respectively, basic vector fields). Indeed let
$X\in\Gamma({\mathcal E}^{4l})$ be a basic vector field. Then we
have immediately
$i_{V}((d\eta_\delta)^{\flat}(X))=d\eta_\delta(X,V)=0$ for any
$V\in\Gamma({\mathcal E}^{4m})$. Next, we have to compute
$i_{V}(d((d\eta_\delta)^{\flat}(X)))(Y)$ for all
$Y\in\Gamma({\mathcal E}^{4l})$. It is not restrictive to assume $Y$
basic. Moreover for simplify the notation we put
$\omega:=(d\eta_\delta)^{\flat}(X)$. Then we have
\begin{align*}
i_{V}(d((d\eta_\delta)^{\flat}(X)))(Y)&=2d\omega(V,Y)=V(\omega(Y))-Y(\omega(V))-\omega([V,Y])\\
&=V(d\eta_\delta(X,Y))-Y(d\eta_\delta(X,V))-d\eta_\delta(X,[V,Y])\\
&=V(d\eta_\delta(X,Y))-Y(d\eta_\delta(X,V))-d\eta_\delta(X,[V,Y])\\
&\quad-X(d\eta_\delta(V,Y))+d\eta_\delta(Y,[V,X])+d\eta_\delta(V,[X,Y])\\
&=3d^2\eta_\delta(X,Y,V)=0,
\end{align*}
for all $V\in\Gamma({\mathcal E}^{4m})$, so that the $1$-form
$(d\eta_\delta)^{\flat}(X)$ is basic. Conversely, let $\omega$ be
a basic $1$-form which vanishes on ${\mathcal E}^{4m+3}$. Then we
prove that, for each $\alpha\in\left\{1,2,3\right\}$, the vector
field $X=(d\eta_\alpha)^{\sharp}(\omega)$ is basic, that is
$[X,V]\in\Gamma({\mathcal E}^{4m})$ for any $V\in\Gamma({\mathcal
E}^{4m})$. Since, by Lemma \ref{lemmapro3},
$[X,V]\in\Gamma({\mathcal H})$, the last condition is equivalent
to require that $d\eta_\alpha([X,V],Y)=0$ for any
$Y\in\Gamma({\mathcal E}^{4l})$. Without loss in generality we can
assume $Y$ to be a basic vector field. We have
\begin{align*}
d\eta_\alpha([X,V],Y)&=3d^2\eta_{\alpha}(V,X,Y)-V(d\eta_\alpha(X,Y))+d\eta_\alpha(X,[V,Y])\\
&\quad-X(d\eta_\alpha(Y,V))-Y(d\eta_\alpha(V,X))+d\eta_\alpha([X,Y],V)\\
&=-V(d\eta_\alpha(X,Y))=-V((d\eta_\alpha)^{\flat}(X)(Y))=-V(\omega(Y))\\
&=-V(\omega(Y))+Y(\omega(V))+\omega([V,Y])=-2d\omega(V,Y)\\
&=-(i_{V}d\omega)(Y)=0
\end{align*}
since $\omega$ is basic.  This proves that $X$ is basic. Therefore
by \eqref{formulaimportante} we get the projectability of
$\phi_\alpha$. Finally we show that ${\mathcal E}^{4m}$ is a
Riemannian foliation, that is for any $V\in\Gamma({\mathcal
E}^{4m})$ $({\mathcal L}_{V}g)|_{N({\mathcal E}^{4m})}=0$, where
$N({\mathcal E}^{4m})=TM/{{\mathcal E}^{4m}}$ is the normal bundle
of the foliation ${\mathcal E}^{4m}$ which is identified with
${\mathcal E}^{4l}\oplus{\mathcal V}$ via the Riemannian metric
$g$. For any $V\in\Gamma({\mathcal E}^{4m})$ and
$X,Y\in\Gamma({\mathcal E}^{4l})$
\begin{align*}
({\mathcal L}_{V}g)(X,Y)&=V(g(X,Y))-g([V,X]_{{\mathcal E}^{4l}},Y)-g(X,[V,Y]_{{\mathcal E}^{4l}})\\
&=-V(d\eta_{\alpha}(X,\phi_\alpha
Y))+d\eta_\alpha([V,X],\phi_\alpha
Y)+d\eta_\alpha(X,\phi_\alpha[V,Y])\\
&=-V(d\eta_{\alpha}(X,\phi_\alpha
Y))+d\eta_\alpha([V,X],\phi_\alpha
Y)+d\eta_\alpha(X,[V,\phi_\alpha Y])\\
&=-({\mathcal L}_{V}d\eta_\alpha)(X,\phi_\alpha Y)=0,
\end{align*}
where we have used the projectability of $d\eta_\alpha$ and
$\phi_\alpha$. Moreover, by Lemma \ref{lemmapro3} and Lemma
\ref{lemmapro4} we get that $({\mathcal
L}_{V}g)(\xi_\delta,Y)=({\mathcal L}_{V}g)(\xi_\delta,\xi_\rho)=0$.
Thus the situation is the following: for each
$\alpha\in\left\{1,2,3\right\}$  $\eta_\alpha$ and $d\eta_\alpha$
project to a $1$-form $\eta'_\alpha$ and a $2$-form
$\Phi'_\alpha=d\eta'_\alpha$; the vector field $\xi_\alpha$ projects
to a vector field $\xi'_\alpha$ satisfying
$\eta'_\alpha(\xi'_\alpha)=1$ and
$d\eta'_\alpha(\xi'_\alpha,\cdot)=0$; the tensor field $\phi_\alpha$
projects to a tensor field $\phi'_\alpha$ such that
$\phi_{\alpha}'^2=-I+\eta'_\alpha\otimes\xi'_\alpha$. Moreover the
Riemannian metric $g$ projects to a Riemannian metric $g'$
compatible with each almost contact structure
$(\phi'_\alpha,\xi'_\alpha,\eta'_\alpha)$. Then one easily  checks
that \eqref{3-sasaki} hold. Finally, that this projected structure
is in fact $3$-$\alpha$-Sasakian follows directly from Corollary
\ref{maximal}.
\end{proof}


We consider now the distribution ${\mathcal E}^{4m+3}={\mathcal E}^{4m}
\oplus\mathcal V$.

\begin{theorem}\label{profinale}
Every $3$-quasi-Sasakian manifold $M^{4n+3}$ of rank $4l+3$ admits a
canonical transversal quaternionic-K\"{a}hler structure given by a
foliation whose leaves are $3$-quasi-Sasakian manifolds of rank $3$.
\end{theorem}
\begin{proof}
By the integrability of ${\mathcal E}^{4m}$ and of $\mathcal V$ and
Lemma \ref{lemmapro4}, it follows that the distribution ${\mathcal
E}^{4m+3}$ is involutive, hence it defines a $(4m+3)$-dimensional
foliation of $M$. Let $N$ be a leaf of this foliation and
$(\phi^N_\alpha,\xi^N_\alpha,\eta^N_\alpha,g^N)$ be the normal
almost $3$-contact metric structure on $N$ obtained from $M$ by
restriction. Then, since $d\Phi_\alpha=0$ and each $1$-form
$\eta_\alpha$ is closed on ${\mathcal E}^{4m}$, we have that $N$ is
endowed with a $3$-quasi-Sasakian structure of rank $3$ (cf.
Proposition \ref{proiettabilita0}). Next, that ${\mathcal E}^{4m+3}$
is a Riemannian foliation follows from the fact that each Reeb
vector field is Killing and ${\mathcal L}_{V}g|_{{\mathcal
E}^{4l}}=0$, $V\in\Gamma({\mathcal E}^{4m})$, (cf. Theorem
\ref{integrabilita}). Finally, let $Q$ be the subbundle of the
endomorphism bundle $\textrm{End}({\mathcal E}^{4l})$ spanned by
$\phi_1$, $\phi_2$, $\phi_3$. Then, since each $\phi_\alpha$ is
foliate with respect to the foliation ${\mathcal E}^{4m}$ (cf.
Theorem \ref{integrabilita}), by \eqref{proiezionephi}, we have that
the subbundle $Q$ is projectable with respect to the foliation
${\mathcal E}^{4m+3}$. Arguing as in \cite{ishihara} one can prove
that the space of leaves is in fact quaternionic-K\"{a}hler.
\end{proof}

Finally, we turn our attention to the distribution ${\mathcal E}^{4l+3}$
on which the 1-forms $\eta_\alpha$ have maximal rank.
\begin{theorem}\label{profinalissimo}
Every $3$-quasi-Sasakian manifold of rank $4l+3$ admits a canonical
transversal hyper-K\"{a}hler structure given by a foliation whose
leaves are $3$-$\alpha$-Sasakian manifolds.
\end{theorem}
\begin{proof}
We first prove that the distribution ${\mathcal E}^{4l+3}$ is
integrable and  defines a Riemannian foliation of the
$3$-quasi-Sasakian manifold
$(M^{4n+3},\phi_\alpha,\xi_\alpha,\eta_\alpha,g)$. Let
$Y,Y'\in\Gamma({\mathcal E}^{4l})$. Then for any
$X\in\Gamma({\mathcal E}^{4m})$ we have
\begin{align*}
0&=3d\Phi_\alpha(X,Y,Y')\\
&=X(\Phi_\alpha(Y,Y'))+Y(\Phi_\alpha(Y',X))+Y'(\Phi_\alpha(X,Y))-\Phi_\alpha([X,Y],Y')\\
&\quad-\Phi_\alpha([Y,Y'],X)-\Phi_\alpha([Y',X],Y)\\
&=({\mathcal L}_{X}\Phi_\alpha)(Y,Y')-\Phi_\alpha([Y,Y'],X)\\
&=({\mathcal L}_{X}g)(Y,\phi_\alpha Y')+g(Y,({\mathcal L}_{X}\phi_\alpha)Y')-\Phi_\alpha([Y,Y'],X)\\
&=-g([Y,Y'],\phi_\alpha X),
\end{align*}
where we have used the projectability of the metric $g$ and of the
tensor field $\phi_\alpha$ with respect to the foliation ${\mathcal
E}^{4m}$, proved in Theorem \ref{integrabilita}. It follows that
$[Y,Y']$ is orthogonal to ${\mathcal E}^{4m}$ and hence belongs to
${\mathcal E}^{4l+3}$. Moreover, by Lemma \ref{lemmapro4} and the integrability of $\mathcal V$, we have that ${\mathcal
E}^{4l+3}$ is integrable and it is easy to check, using
\eqref{formulapesante}, that the normal almost $3$-contact metric
structure induced from $(\phi_\alpha,\xi_\alpha,\eta_\alpha,g)$ on
each leaf of ${\mathcal E}^{4l+3}$ is in fact $3$-$\alpha$-Sasakian.
Now we pass to study the space of leaves of ${\mathcal E}^{4l+3}$.
We prove that ${\mathcal E}^{4l+3}$ is a Riemannian foliation and
that the tensor fields $\phi_1$, $\phi_2$, $\phi_3$ locally project,
together with $g$, to a hyper-K\"{a}hler structure on the space of
leaves. For all $Y\in\Gamma({\mathcal E}^{4l+3})$ and
$X,X'\in\Gamma({\mathcal E}^{4m})$ we have
\begin{align*}
({\mathcal
L}_{Y}\Phi_\alpha)(X,X')&=Y(\Phi_\alpha(X,X'))-\Phi_\alpha([Y,X],X')-\Phi_\alpha(X,[Y,X'])\\
&=X(\Phi_\alpha(X',Y))+X'(\Phi_\alpha(Y,X))+Y(\Phi_\alpha(X,X'))\\
&\quad-\Phi_\alpha([X,X'],Y)-\Phi_\alpha([X',Y],X)-\Phi_\alpha([Y,X],X')\\
&=3d\Phi_\alpha(X,X',Y)=0,
\end{align*}
thus each fundamental $2$-form $\Phi_\alpha$ projects to a  $2$-form
$\Omega'_\alpha$ which is closed since $\Phi_\alpha$ is. Next we
prove that also each tensor field $\phi_\alpha$ is foliate, that is
it maps basic vector fields to basic vector fields. As in Lemma 4.1
in \cite{cappellettidenicola} we have that, for an even permutation
$(\alpha,\beta,\gamma)$ of $\left\{1,2,3\right\}$,
\begin{equation}\label{formulaimportante1}
\phi_\alpha X=-(\Phi_\beta)^{\sharp}(\Phi_\gamma)^{\flat}X
\end{equation}
for all $X\in\Gamma({\mathcal E}^{4m})$, where, for each
$\delta\in\left\{1,2,3\right\}$, $(\Phi_\delta)^{\flat}:X \mapsto
\Phi_\delta(X,\cdot)$ is the musical isomorphism induced from
$\Phi_\delta$ between ${\mathcal E}^{4m}$ and the $1$-forms which
vanish on ${\mathcal E}^{4l+3}$, and $(\Phi_\delta)^{\sharp}$
denotes its inverse. Therefore, in order to prove that $\phi_\alpha$
is foliate it is sufficient to check that, for each
$\delta\in\left\{1,2,3\right\}$, $(\Phi_\delta)^{\flat}$ maps basic
vector fields to basic $1$-forms and, conversely,
$(\Phi_\delta)^{\sharp}$ maps basic $1$-forms to basic vector
fields. Let $X\in\Gamma({\mathcal E}^{4m})$ be a basic vector field.
We have to show that the $1$-form $\omega:=(\Phi_\delta)^{\flat}X$
is basic, i.e. satisfies $i_{Y}\omega=i_{Y}d\omega=0$ for all
$Y\in\Gamma({\mathcal E}^{4l+3})$. Indeed we have
$i_{Y}\omega=\omega(Y)=\Phi_\delta(X,Y)=g(X,\phi_\delta Y)=0$ since
$\phi_\delta({\mathcal E}^{4l+3})\subset{\mathcal E}^{4l+3}$. Next,
one has
$i_{Y}d\omega(X')=2d\omega(Y,X')=Y(\omega(X'))-X'(\omega(Y))-\omega([Y,X'])=({\mathcal
L}_{Y}\Phi_\delta)(X,X')=0$ for any $X'\in\Gamma({\mathcal E}^{4m})$
(which is not restrictive to assume basic). Conversely, for any
basic $1$-form $\omega$ we have to show that the vector field
$X:=(\Phi_\delta)^{\sharp}(\omega)$ is basic, that is
$[X,Y]\in\Gamma({\mathcal E}^{4l+3})$ for any $Y\in\Gamma({\mathcal
E}^{4l+3})$. This last condition is equivalent to require that
$\Phi_\delta([X,Y],X')=0$ for any $X'\in\Gamma({\mathcal E}^{4m})$.
It is not restrictive to assume $X'$ basic. Then we have
\begin{align*}
\Phi_\delta([X,Y],X')&=-3d\Phi_\delta(X,Y,X')+X(\Phi_\delta(Y,X'))+Y(\Phi_\delta(X',X))\\
&\quad+X'(\Phi_\delta(X,Y))-\Phi_\delta([Y,X'],X)-\Phi_\delta([X',X],Y)\\
&=-Y(\Phi_\delta(X,X'))\\
&=i_{Y}d\omega(X')=0.
\end{align*}
It remains to prove that the Riemannian metric $g$ is bundle-like.
This follows easily from the projectability of $\Phi_\alpha$ and
$\phi_\alpha$. Indeed for any $Y\in\Gamma({\mathcal E^{4l}})$ and
$X,X'\in\Gamma({\mathcal E^{4m}})$ we have
\begin{equation*}
({\mathcal L}_{Y}g)(X,X')=-({\mathcal
L}_{Y}\Phi_\alpha)(X,\phi_\alpha X')+g(X,({\mathcal
L}_{Y}\phi_\alpha)\phi_\alpha X')=0,
\end{equation*}
whereas $({\mathcal L}_{\xi_\alpha}g)(X,X')=0$ since $\xi_\alpha$ is
Killing. We denote by $J'_\alpha$ and $g'$ the tensor fields induced on the space of leaves by each
$\phi_\alpha$ and by the Riemannian metric $g$. Then a
straightforward computation yields that
$(J'_\alpha,\Omega'_\alpha,g')$ is an almost hyper-Hermitian
structure. Thus, the closedness of $\Omega'_1$, $\Omega'_2$,
$\Omega'_3$ imply, by the Hitchin Lemma (\cite{hitchin}), that
$(J'_\alpha,\Omega'_\alpha,g')$ is in fact hyper-K\"{a}hler.
\end{proof}

\begin{corollary}\label{prodotto}
Let $(M^{4n+3},\phi_\alpha,\xi_\alpha,\eta_\alpha,g)$ be a
$3$-quasi-Sasakian manifold of rank $4l+3$ with
$\left[\xi_\alpha,\xi_\beta\right]=c\xi_\gamma$, $c\neq 0$. Then
$M^{4n+3}$ is locally the Riemannian product of a
$3$-$\alpha$-Sasakian manifold $M^{4l+3}$, where
$\alpha=\frac{c}{2}$, and a hyper-K\"{a}hler manifold $M^{4m}$, with
$m=n-l$.
\end{corollary}
\begin{proof}
The tangent bundle of $M^{4n+3}$ splits up as the orthogonal sum of
the Riemannian foliations ${\mathcal E}^{4l+3}$ and ${\mathcal
E}^{4m}$. Because of the duality Riemannian-totally geodesic,
${\mathcal E}^{4l+3}$ and ${\mathcal E}^{4m}$ are also totally
geodesic foliations. It follows that $M^{4n+3}$ is the Riemannian
product of a leaf $M^{4l+3}$ of ${\mathcal E}^{4l+3}$ and a leaf
$M^{4m}$ of ${\mathcal E}^{4m}$. Taking into account that
$\psi_\alpha$ and $\phi_\alpha$ agree on ${\mathcal E}^{4l+3}$ and
applying \eqref{formulapesante}, we have that
$(\psi_\alpha,\xi_\alpha,\eta_\alpha,{g})|_{{\mathcal E}^{4l+3}}$ is
an almost $3$-$\alpha$-Sasakian structure over $M^{4l+3}$, where we
have put $\alpha=\frac{c}{2}$. Hence, by Proposition
\ref{lemmamino}, it is $3$-$\alpha$-Sasakian. Since $\theta_\alpha$
agrees with $\phi_\alpha$ on ${\mathcal E}^{4m}$, the maps
$\theta_\alpha|_{{\mathcal E}^{4m}}$ define a quaternionic structure
which is compatible with the metric ${g}|_{{\mathcal E}^{4m}}$.
Finally, define the $2$-forms ${\Theta}_\alpha$ by
${\Theta}_\alpha\left(X,Y\right)={g}\left(X,\theta_\alpha Y\right)$
for any $X,Y\in\Gamma({{\mathcal E}^{4m}})$. We have
${\Theta}_\alpha={\Phi}_\alpha|_{{\mathcal E}^{4m}}$ and hence
$d{\Theta}_\alpha=0$.  By virtue of the mentioned Hitchin Lemma
(\cite{hitchin}) the structure defined on $M^{4m}$ turns out to be
hyper-K\"{a}hler.
\end{proof}

\begin{remark}
Note that Corollary \ref{prodotto} strongly improves, both in the
assumptions and in the results, the splitting theorem which is
proven in \cite{mag} (Theorem 5.6). It should be also emphasized
that an analogous result does not hold for a quasi-Sasakian
manifold.
\end{remark}

Corollary \ref{prodotto} also happens to have some notable
consequences. First of all, the structure group of a
$3$-quasi-Sasakian manifold of rank $4l+3$ is clearly reduced to
$\emph{Sp}(m)\times\emph{Sp}(l) \times I_3$. Next, we obtain an
improving of Theorem \ref{pro1}. Namely, under the assumption of
regularity for the foliation $\mathcal V$, the space of leaves
$M/{\mathcal V}$ is an almost quaternionic-Hermitian manifold which
is the local Riemannian product of a quaternionic-K\"{a}hler
manifold and a hyper-K\"{a}hler manifold. Another consequence is
related to the concept of contact-symplectic pair defined by Bande
in \cite{bande2003}. A \emph{contact-symplectic pair} of type $(h,k)$ on a
manifold $M$ of dimension $2h+2k+1$ is a pair $(\beta,\omega)$ of a
$1$-form $\beta$ and a closed $2$-form $\omega$ such that the
following conditions are satisfied: $\beta \wedge (d \beta)^h \wedge
\omega^k$ is a volume form, $(d \beta)^{h+1}=0$ and
$\omega^{k+1}=0$. Due to Corollary \ref{prodotto}, any
$3$-quasi-Sasakian manifold of rank $4l+3$ is canonically endowed
with nine contact-symplectic pairs as it is stated in the following
claim which is easy to verify.
\begin{corollary}\label{coppie}
Let $(M^{4n+3},\phi_\alpha,\xi_\alpha,\eta_\alpha,g)$ be a
$3$-quasi-Sasakian manifold of rank $4l+3$. Then, for each $\beta, \gamma\in\{1,2,3\}$,
$(\eta_\beta,\Theta_{\gamma})$ is a contact-symplectic pair of type $(2l+1,2m)$.
\end{corollary}
%
%
Finally, using Corollary \ref{prodotto} we can compute the complete
expression of the Ricci tensor in any $3$-quasi-Sasakian manifold.
Before, we prove the following preliminary result.

\begin{proposition}\label{eistein1}
Every $3$-$\alpha$-Sasakian manifold of dimension $4n+3$ is an
Einstein manifold with Einstein constant $2\alpha^2(2n+1)$.
\end{proposition}
\begin{proof}
Let $(M,\phi_\delta,\xi_\delta,\eta_\delta,g)$,
$\delta\in\left\{1,2,3\right\}$, be a $3$-$\alpha$-Sasakian
manifold. Then by virtue of \eqref{conseguenze} and Proposition 4.4
in \cite{olszak2} we have that
$(\phi_\delta,\xi_\delta,\eta_\delta,g)$ can be obtained by a
homothetic deformation of a $3$-Sasakian structure
$(\bar{\phi}_\delta,\bar{\xi}_\delta,\bar{\eta}_\delta,\bar{g})$
given by
\begin{equation*}
\bar\phi_\delta=\phi_\delta, \
\bar\xi_\delta=\frac{1}{\alpha}\xi_\delta, \
\bar\eta_\delta=\alpha\eta_\delta, \ \bar g=\alpha^2 g.
\end{equation*}
Then, since it is well known that any $3$-Sasakian manifold is
Einstein, we conclude that also the metric $g$ is Einstein. For
computing the Einstein constant $\lambda$ we use \eqref{olszakricci}
and \eqref{conseguenze} getting
\begin{equation*}
\lambda=\lambda
g(\xi_\delta,\xi_\delta)=\textrm{Ric}(\xi_\delta,\xi_\delta)=\|\nabla\xi_\delta\|^2=2\alpha^2(2n+1).
\end{equation*}
\end{proof}

As pointed out in Corollary \ref{riccicurv1}, $3$-quasi-Sasakian manifolds of rank $4l+1$ are Ricci-flat, being
$3$-cosymplectic. As for the case of rank $4l+3$, we have the following

\begin{theorem}\label{eistein2}
Let $M$ be a $3$-quasi-Sasakian manifold of rank $4l+3$. The Ricci
tensor is given by
\begin{equation}\label{espressionericci1}
\emph{Ric}(X,Y)= \left\{
  \begin{array}{ll}
    \frac{c^2}{2}(2l+1)g(X,Y), & \hbox{if $X,Y\in\Gamma({\mathcal E}^{4l+3})$;} \\
    0, & \hbox{elsewhere.}
  \end{array}
\right.
\end{equation}
In particular, $M$ has positive scalar
curvature $\frac{c^2}{2}(2l+1)(4l+3)$.
\end{theorem}
\begin{proof}
In view of Corollary \ref{prodotto}, $M$ is locally the Riemannian
product of a $(4l+3)$-dimensional $3$-$\alpha$-Sasakian manifold
$M'$, with $\alpha=\frac{c}{2}$, and a $4m$-dimensional
hyper-K\"{a}hler manifold $M''$. Thus, because of Proposition
\ref{eistein1} and the Ricci-flatness of $M''$, we get the
assertion.
\end{proof}



\begin{corollary}
No $3$-quasi-Sasakian manifold is $\eta$-Einstein unless the
following cases:
\begin{enumerate}
  \item[(i)] $3$-$\alpha$-Sasakian manifolds, which are Einstein with
strictly positive scalar curvature;
  \item[(ii)] $3$-cosymplectic manifolds, which are Ricci-flat;
  \item[(iii)] $3$-quasi-Sasakian manifolds of rank $3$, which are
$\eta$-Einstein non-Einstein.
\end{enumerate}
\end{corollary}

\begin{remark}
Applying the Pasternack's refinement \cite{pasternack} of the
classical Bott vanishing theorem to the Riemannian foliations
$\mathcal V$, ${\mathcal E}^{4m}$, ${\mathcal E}^{4m+3}$ and
${\mathcal E}^{4l+3}$, considered in Theorem \ref{pro1}, Theorem
\ref{integrabilita}, Theorem \ref{profinale}, and Theorem
\ref{profinalissimo}, respectively, we get the following topological
obstructions to the existence of a $3$-quasi-Sasakian structure of
rank $4l+3$ on a manifold $M$ of dimension $4n+3$:
$\textrm{Pont}^j({\mathcal H})=0$ for all $j>4n$,
$\textrm{Pont}^j({\mathcal E}^{4l+3})=0$ for all $j>4m$,
$\textrm{Pont}^j({\mathcal E}^{4l})=0$ for all $j>4m+3$ and
$\textrm{Pont}^j({\mathcal E}^{4m})=0$ for all $j>4l+3$ , where
$\textrm{Pont}({\mathcal H})$, $\textrm{Pont}({\mathcal
E}^{4l+3})=0$, $\textrm{Pont}({\mathcal E}^{4l})$ and
$\textrm{Pont}({\mathcal E}^{4m})$ denote, respectively, the
Pontryagin algebras of the subbundles $\mathcal H$, ${\mathcal
E}^{4l+3}$, ${\mathcal E}^{4l}$ and ${\mathcal E}^{4m}$ of the
tangent bundle of $M$. Furthermore, the vanishing of these primary
characteristic classes permits also the construction of some
secondary characteristic classes as it is done in \cite{lazarov} and
\cite{piccinnivaisman}.
\end{remark}

\section*{Acknowledgments}
The second author acknowledges financial support by a CMUC
postdoctoral fellowship.

\end{document}